\documentclass[11pt]{article}
\usepackage[a4paper]{geometry}
\usepackage{amsthm}
\usepackage{amsmath}
\usepackage{amssymb}
\usepackage{amsfonts}
\usepackage{epsfig,graphicx}

\begin{document}

\theoremstyle{theorem}
\newtheorem{Lemma}{Lemma}[section]
\newtheorem{Proposition}{Proposition}[section]
\newtheorem{Corollary}{Corollary}[section]
\newtheorem{Theorem}{Theorem}[section]
\newtheorem{Condition}{Condition}[section]
\newtheorem{Claim}{Claim}[section]
\newtheorem{Integrator}{Integrator}

\theoremstyle{definition}
\newtheorem{Definition}{Definition}[section]
\newtheorem{Remark}{Remark}[section]
\newtheorem{Example}{Example}[section]

\title{Space-time FLAVORS: finite difference, multisymlectic, and pseudospectral integrators for multiscale PDEs}

\author{Molei Tao, Houman Owhadi, Jerrold E. Marsden}

\maketitle

\begin{abstract}
We present a new class of integrators for stiff PDEs. These integrators are generalizations of  FLow AVeraging integratORS (FLAVORS) for stiff ODEs and SDEs introduced in \cite{FLAVOR10} with the following properties: (i) {\it Multiscale}: they are based on flow averaging and have a computational cost determined by mesoscopic steps in space and time instead of microscopic steps in space and time; (ii) {\it Versatile}: the method is based on averaging the flows of the given PDEs (which may have hidden slow and fast processes). This bypasses the need for identifying explicitly (or numerically) the slow variables or reduced effective PDEs; (iii) {\it Nonintrusive}: A pre-existing numerical scheme resolving the microscopic time scale can be used as a black box and easily turned into one of the integrators in this paper by turning the large coefficients on over a microscopic timescale and off during a mesoscopic timescale; (iv) {\it Convergent over two scales}: strongly over slow processes and in the sense of measures over fast ones; (v) {\it Structure-preserving}: for stiff Hamiltonian PDEs (possibly on manifolds), they can be made to be multi-symplectic,  symmetry-preserving (symmetries are group actions that leave the system invariant) in all variables and variational.
\end{abstract}

\section{Introduction}
Multi-scale PDEs can be divided into two (possibly over-lapping) categories: PDEs with highly oscillating or rough coefficients and PDEs with large (or stiff) coefficients. Classical numerical methods  are usually: (i) stable but arbitrarily inaccurate for the former category (consider, for instance, a finite element method for the elliptic operator $-\operatorname{div}(a\nabla)$ with a rapidly changing coefficient $a\in L^\infty$), or (ii) unstable for the latter category. Accurate numerical methods for the former category, called numerical homogenization methods, are, in absence of local ergodicity or scale separation, based on the compactness of the solution space (we refer, for instance, to \cite{MR2292954, BerlyandOwhadi10, OwZhLoc10}).
Numerical methods for the latter category are, in essence, based on the existence of slow and fast variables (or components) \cite{MR2164093}. When fast variables converge toward Dirac (single point) distributions, asymptotic-preserving schemes \cite{FiJi10} allow for simulations with large time steps.
We also refer to \cite{MR1781209, MR1756424} for multi-scale transport equations and hyperbolic systems of conservation laws with stiff diffusive relaxation.
 Well-identified slow variables can be simulated with large time-steps using the two-scale structure of the original stiff PDEs (we refer to \cite{ArGeKeSlTi09} and \cite{MR2314852} for existing examples; slow variables satisfy a non-stiff PDE that can be identified in analogy to equations (A.9) and (A.13) of \cite{FLAVOR10}; we also refer to \cite{MR2164093} for a definition of slow variables).

In this paper, we consider the second category of PDEs and propose a generalization of FLow AVeraging integratORS (FLAVORS) (introduced in \cite{FLAVOR10} for stiff ODEs and SDEs) to stiff PDEs. Multi-scale integrators for stiff PDEs are obtained without the identification of slow variables by turning on and off stiff coefficients in single-step (legacy) integrators (used as black boxes) and alternating microscopic and mesoscopic time steps (Subsection \ref{generalmeth}). We illustrate the generality of the proposed strategy by applying it to  finite difference methods in Section \ref{FDsection}, multi-symplectic integrators in Section \ref{multi-symplecticSection}, and pseudospectral methods in Section \ref{secpseudospectral} (although we have not done so in this paper, the proposed strategy can also be applied to finite element methods or finite volume methods).
The convergence of the proposed strategy, after semi-discretization in space, is analyzed in Subsection \ref{jhsjgfdhgdfdhgf3}, where a non-asymptotic error bound indicates the two-scale convergence (\cite{FLAVOR10}, i.e., strong with respect to hidden slow variables and weak with respect to hidden fast variables) of PDE-FLAVORS. As illustrated by numerical (Figure \ref{GLerror}) and theoretical results (Section \ref{SectionErrorAnalysis}), an explicit tuning ($(h/\epsilon)^2 \ll H \ll h/\epsilon$) between microscopic $h$ and mesoscopic ($H$) time-steps and the stiff parameter $1/\epsilon$ is necessary and sufficient for convergence.
We also show in  Section \ref{seccharact} that applying the FLAVOR strategy to characteristics leads to accurate approximations of solutions of stiff PDEs.

These results, along with those of \cite{FLAVOR10}, diverge from the concept that, in situations where the slow variables are not linear functions of the original variables, multiscale algorithms \emph{``do not work'' ``if the slow variables are not explicitly identified and made use of''} (page 2 of \cite{Nested07}).

\section{Finite difference and space-time FLAVOR mesh}
\label{FDsection}

\subsection{Single-scale method and limitation}
Consider a multiscale PDE:
\begin{equation}
    F(1, \epsilon^{-1}, x,t,u(x,t),u_x(x,t),u_t(x,t),u_{xx}(x,t),u_{xt}(x,t),u_{tt}(x,t), \ldots)=0
    \label{generalPDE}
\end{equation}
where $F$ is a given function (possibly nonlinear), $\epsilon$ is a small positive real parameter and $x$ and $t$ are spatial and temporal coordinates.

 To obtain a numerical solution of \eqref{generalPDE}, the simplest single-scale finite difference approach employs a uniform rectangular mesh with time step length $h$ and space step length $k$, and approximates the solution $u$ at its values at discrete grid points. Differential operators will be approximated by finite differences; for instance, according to forward space forward time rules: $u_x(ik,jh) \approx (u_{i+1,j}-u_{i,j})/k$ and $u_t(ik,jh) \approx (u_{i,j+1}-u_{i,j})/h$, where $u_{ij}$ is the numerical solution at discrete grid point with space index $i$ and time index $j$. After this discretization, the original PDE is approximated by a finite dimensional algebraic system, which can be solved to yield the numerical solution.

Of course, a necessary condition for obtaining stability and accuracy in the numerical solution is that $h$ and $k$ have to be small enough. A quantitative statement on how small they need to be will depend on the specific PDE and discretization. For 1D linear advection equations $u_x-au_t=0$ and forward time forward space discretizations, the $h<k/a$  CFL condition \cite{CFL} has to be met to ensure stability, which is also a neccessary condition for accuracy \cite{LaRi56}. Intuitively, the CFL condition guarantees that information does not propagate faster than what the numerical integrator can handle. The Von Neumann stability analysis \cite{ChFjNe50} helps determine analogous CFL conditions for linear equations with arbitrary discretizations. The stability of numerical schemes for general nonlinear equations remains a topic of study. We refer to \cite{St04} for additional discussions on single-scale finite difference schemes.  In general, the presence of a stiff coefficient $\epsilon^{-1}$ in equation \eqref{generalPDE} requires  $h$ and $k$  to scale with $\epsilon$ in order to guarantee the stability of numerical integration schemes.
 This makes the numerical approximation of the solution of \eqref{generalPDE} computationally untractable
  when $\epsilon$ is close to 0.

\subsection{Multiscale FLAVORization and general methodology}\label{generalmeth}

FLAVORs are multiscale in the sense that they accelerate computation by adopting both larger time and space steps. A finite difference scheme can be FLAVORized by employing two rules:

\begin{figure} [ht]
\vspace{-10pt}
\includegraphics[width=\textwidth]{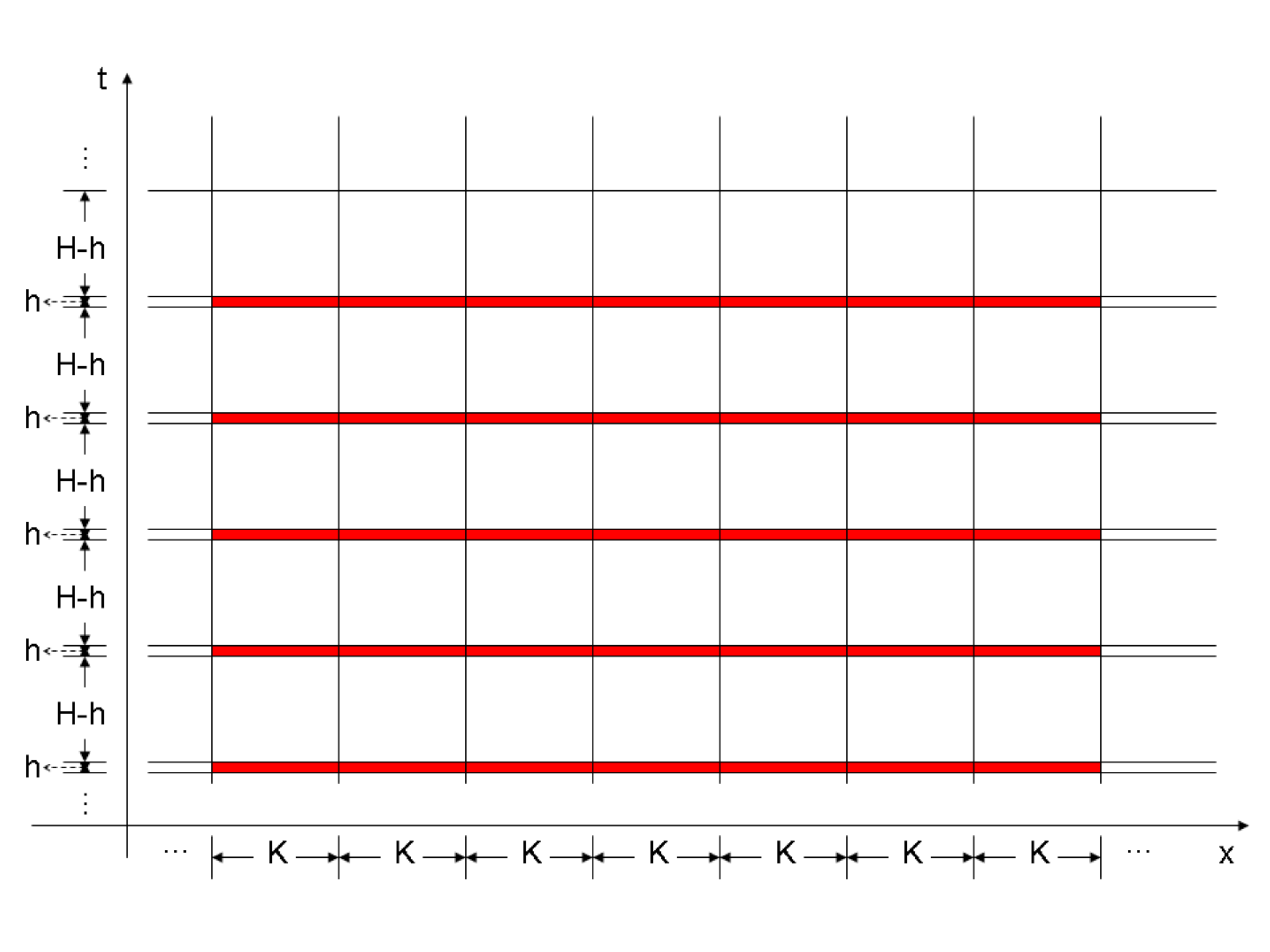}
\vspace{-40pt}
\caption{\footnotesize Mesh used by FLAVORS. A uniform mesoscopic space step is used and two alternating small and mesoscopic time steps are used. Stiffness is turned on in red regions and turned off otherwise.}
\label{FLAVORmesh}
\end{figure}

First, instead of a uniform mesh, use a mesh as depicted in Figure \ref{FLAVORmesh}, in which a uniform spatial grid corresponds to a mesoscopic space step $K$ that does not scale with $\epsilon$, and an alternating temporal grid corresponds to two time steps, microscopic $h$ (scaling with $\epsilon$) and mesoscopic $H-h$ ($H$ independent from $\epsilon$). It is worth mentioning that when using this non-uniform mesh, grid sizes have to be taken into consideration when derivatives are approximated by finite differences. 1st-order derivatives are straightforward to obtain, and we refer to Section \ref{multi-symplecticSection} for approximations of higher order derivatives.

Second, the stiff parameter $\epsilon^{-1}$ should be temporarily set to be 0 (i.e., turned off) when the current time step is the mesoscopic $H-h$; if the small time step $h$ is used instead, the large value of $\epsilon^{-1}$ needs to be restored, or in other words, stiffness should be turned on again.

The rule of thumb is that $k$ and $h$ should be chosen such that the integration of \eqref{generalPDE} with these step sizes and stiffness turned on is stable and accurate. On the other hand, there is another pair of step size values such that the same integration with stiffness turned off is stable and accurate, and $K$ and $H$ should be chosen to be an order of magnitude smaller than these values. FLAVORS does not require a microscopic $k$, but only a mesoscopic space-step $K$, a microscopic time-step $h$, and a mesoscopic time-step $H$.

The intuition is as follows: adopt the  point of view of semi-discrete approach for PDE integration, in which space is discretized first and the PDE is approximated by a system of ODEs. The integration (in the time) of the resulting finite dimensional ODE system can be accelerated by applying the FLAVOR strategy to any legacy scheme (used as a black box). Turning on and off stiff coefficients in the legacy scheme and alternating microscopic time steps (stiffness on) with mesoscopic time steps (stiffness on) preserves the symmetries of that scheme and at the same time induces an averaging of the dynamic of (possibly hidden) slow variables with respect to the fast ones. With this strategy, the FLAVORized scheme advances in mesoscopic time steps without losing stability. The (possibly hidden) slow dynamic is captured in a strong sense, while the fast one is captured only in the (weak) sense of measures.
A rigorous proof of convergence of the proposed method relies on the assumption of existence of (possibly hidden) slow variables and of local ergodicity of (possibly hidden) fast variables (we refer to Section \ref{SectionErrorAnalysis}). It is important to observe that the proposed method does not require the identification of slow variables.

\subsection{Example: conservation law with Ginzburg-Landau source}
Consider a specific stiff PDE:
\begin{equation}
    u_t+f(u)_x=\epsilon^{-1} u(1-u^2)
    \label{GLeq}
\end{equation}
in which $f(u)=\sin u$ and $0 < \epsilon \ll 1$. Use the boundary condition of $u(x=0,t)=u(x=L,t)$ and the initial condition of $u(x,t=0)=\sin(\pi x)$. This system contains two scales: the fast process corresponds to $u$ quickly converging towards $1$ or $-1$, and the slow process corresponds to the front (with steep gradients) separating $u>0$ from $u<0$ propagating at an $\mathcal{O}(1)$ velocity.

We will FLAVORize the following Lax-Friedrichs finite difference scheme:
\begin{equation}
\begin{cases}
    u_{i+1,j+1} &= \bar{u}_{i+1,j} -h \left( f_u(\bar{u}_{i+1,j}) \frac{u_{i+2,j}-u_{i,j}}{2k} + \epsilon^{-1} \bar{u}_{i+1,j}(1-\bar{u}_{i+1,j}^2) \right) \\
    \bar{u}_{i+1,j} &\triangleq \frac{u_{i+2,j}+u_{i,j}}{2}
    \label{GL_FD}
\end{cases}
\end{equation}
where $u_{i,j}=u_{i+L/k,j}$ and $u_{i,1}=\sin\left(\pi(i-1)k\right)$. If the domain of integration is restricted to $[0,L]\times[0,T]$, then $i=1,2,\ldots,\lfloor L/k \rfloor+1$, and $j=1,2,\ldots,\lfloor T/h \rfloor+1$. We use $h=0.1\epsilon$ and $k=0.2\epsilon$ for our purposes, both of which we found numerically at the order of the stability limit. In our experiment, we chose $\epsilon=2\cdot 10^{-3}$, and therefore $h=0.0002$ and $k=0.0004$.

The FLAVORized version of this scheme is:
\begin{equation}
    \begin{cases}
    \tilde{u}_{i+1,j} &= \bar{u}_{i+1,j}-h \left( f_u(\bar{u}_{i+1,j}) \frac{u_{i+2,j}-u_{i,j}}{2K} + \epsilon^{-1} \bar{u}_{i+1,j}(1-\bar{u}_{i+1,j}^2) \right) \\
    \bar{u}_{i+1,j} &\triangleq (u_{i+2,j}+u_{i,j})/2 \\
    u_{i+1,j+1} &= \frac{\tilde{u}_{i+2,j}+\tilde{u}_{i,j}}{2}-(H-h) \left( f_u(\frac{\tilde{u}_{i+2,j}+\tilde{u}_{i,j}}{2}) \frac{\tilde{u}_{i+2,j}-\tilde{u}_{i,j}}{2K} \right)
    \end{cases}
    \label{GL_FLAVOR}
\end{equation}
where $u_{i,j}=u_{i+L/K,j}$ and $u_{i,1}=\sin\left(\pi(i-1)K\right)$. If the domain of integration is restricted to $[0,L]\times[0,T]$, then $i=1,2,\ldots,\lfloor L/K \rfloor+1$, and $j=1,2,\ldots,\lfloor T/H \rfloor+1$. We use the same $h$ as before, and choose $H=0.005$ and $K=0.01$, which ensures that the stability of the integration remains independent of $\epsilon$.

\begin{figure} [ht]
\begin{tabular}{c}
\hspace{-30pt}
\includegraphics[width=\textwidth]{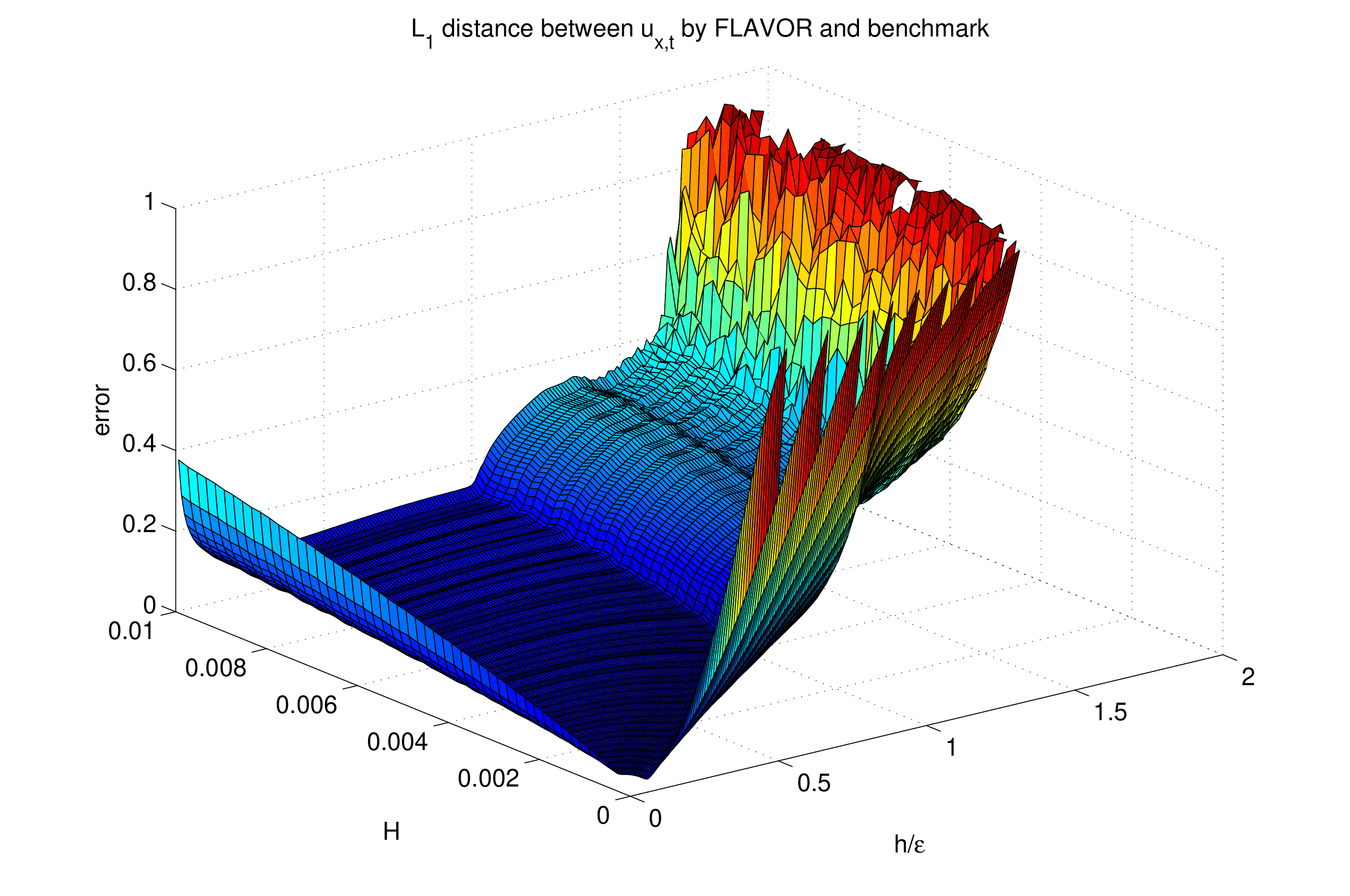}
\hspace{-20pt}
\end{tabular}
\caption{\footnotesize Errors of FLAVOR based on Lax-Friedrichs as a function of $H$ and $h$. $H$ samples multiples of $0.1\epsilon$, starting from 2x to 50x with 1x increment, and $h$ ranges from $0.01\epsilon$ to $3\epsilon$ with $0.01\epsilon$ increment. Errors with magnitude bigger than 1 are not plotted, for they indicate unstable integrations.}
\label{GLerror}
\end{figure}

\begin{figure} [ht]
\begin{tabular}{cc}
\hspace{-30pt}
\includegraphics[width=0.58\textwidth]{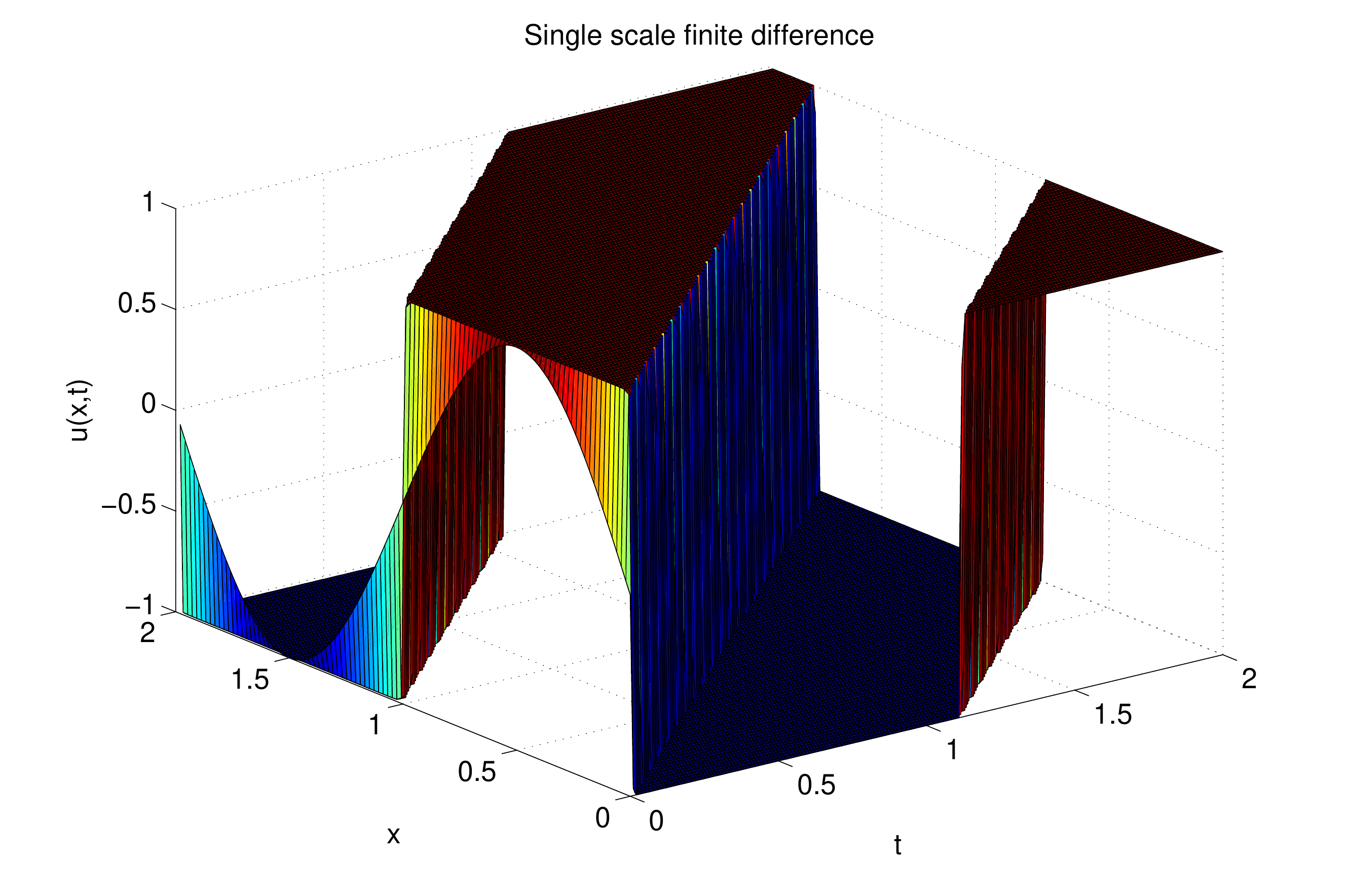}
\hspace{-20pt}
&
\hspace{-20pt}
\includegraphics[width=0.58\textwidth]{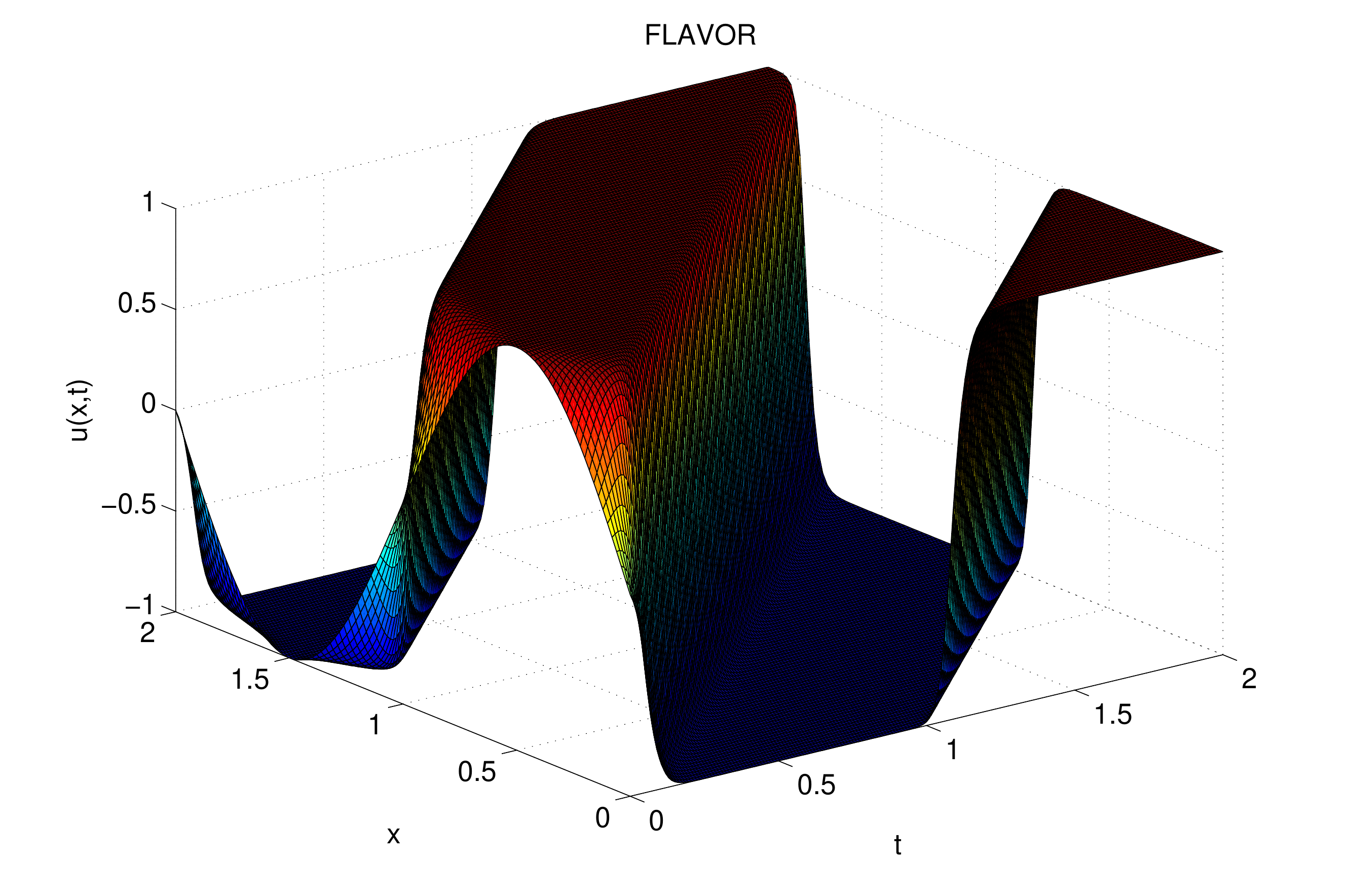}
\hspace{-30pt}
\end{tabular}
\caption{\footnotesize Numerical solutions to \eqref{GLeq} by Lax-Friedrichs (left, Eq. \ref{GL_FD}) and its FLAVORization (right, Eq. \ref{GL_FLAVOR}).}
\label{GLresult}
\end{figure}

Errors of FLAVOR based on Lax-Friedrichs with different $H$ and $h$ values are computed by comparing the results to a benchmark Lax-Friedrichs integration with fine steps $h=0.1\epsilon$ and $k=0.2\epsilon$. More precisely, we calculated the distance between two vectors respectively corresponding to FLAVOR and Lax-Friedrichs integrations, which contain ordered $u(x,t)$ values on the intersection of FLAVOR and Lax-Friedrichs meshes (which is in fact the FLAVOR mesh as long as $H$ is a multiple of $0.1\epsilon$). 1-norm is used and normalized by the number of discrete points to mimic the $L^1$ norm for the continuous solution. Experimental settings are $\epsilon=2\cdot 10^{-3}$, $L=2$ and $T=2$. As we can see in Figure \ref{GLerror}, FLAVOR is indeed uniformly convergent in the sense that the error scales with $H$, as long as $h$ takes an appropriate value. This is not surprising, because we have already proven in the ODE case that the error is bounded by a function of $H$ (uniformly in $\epsilon$) as long as $\left(\frac{h}{\epsilon}\right)^2 \ll H \ll h/\epsilon$, and this error can be made arbitrarily small as $H\downarrow 0$  (notice $H$ can still be much larger than $\epsilon$ as $\epsilon\downarrow 0$).

Also, a typical run of FLAVOR ($H=0.005$ and $K=0.01$) in comparison to the benchmark ($h=0.0002$ and $k=0.0004$) is shown in Figure \ref{GLresult}. FLAVOR captured the slow process strongly in the sense that it obtained the correct speeds of both steep gradients' propagations (up to arithmetic error and fringing).  In this setting, FLAVOR achieves a $\frac{HK}{2hk}=312.5$ fold acceleration. It is worth restating that both spatial and temporal step lengths of FLAVOR are mesocopic, whereas the counterparts in a single scale finite difference method have to be both microscopic for stability. The computational gain by FLAVOR will go to infinity as $\epsilon \rightarrow 0$, and this statement will be true for all FLAVOR examples shown in this paper.

\section{Multisymplectic integrator for Hamiltonian PDEs}
\label{multi-symplecticSection}
\subsection{Single-scale method}
We refer to \cite{BrRe01, MaPaSh98, MaSh99} for a discussion on the geometry of Hamiltonian PDEs (e.g., multi-symplectic structure). We will now recall the Euclidean coordinate form of a Hamiltonian PDE:
\begin{equation}
    \mathcal{M}z_t+\mathcal{K}z_x=\nabla_z H(z)
    \label{HamiltonianPDE}
\end{equation}
where $z(x,t)$ is a n-dimensional vector, $\mathcal{M}$ and $\mathcal{K}$ are arbitrary skew-symmetric matrices on $\mathbb{R}^n$, and $H:\mathbb{R}^n \rightarrow \mathbb{R}$ is an arbitrary smooth function. The solution preserves the multi-symplectic structure in the following sense:
\begin{equation}
    \partial_t \iota(U,V)+ \partial_x \kappa(U,V) = 0
\end{equation}
where $\iota$ and $\kappa$ are differential 2-forms defined by
\begin{equation}
    \iota(x,y)=\langle \mathcal{M}x,y \rangle \quad \text{and} \quad \kappa(x,y)=\langle \mathcal{K}x,y \rangle
\end{equation}
and $U$ and $V$ are two arbitrary solutions to the variational equation (the solution is identified with $dz:\mathbb{R}^2\mapsto \mathbb{R}^n$):
\begin{equation}
    \mathcal{M}dz_t+\mathcal{K}dz_x=D_{zz} H(z) dz,  \quad dz(x,t) \in \mathbb{R}^n
\end{equation}
Preservation of multi-symplecticity can be partially and intuitively interpreted as a conservation of infinitesimal volume in the jet bundle, which generalizes the conservation of phase space volume in Hamiltonian ODE settings to field theories.

A broad spectrum of PDEs fall in the class of Hamiltonian PDEs, including generalized KdV, nonlinear Schr\"{o}dinger models, nonlinear wave equations, atmospheric flows, fluid-structure interactions, etc. \cite{Br97a,Br97b,BrDe99,BrRe01}. We also refer to \cite{BrRe06} and references therein for surveys on numerical recipes, and to \cite{LeMaOrWe2003} for an application to numerical nonlinear elastodynamics.

Hamiltonian PDEs \eqref{HamiltonianPDE} can be viewed as Euler-Lagrange equations for field theories, which are obtained by applying Hamilton's principle (i.e., a variational principle of $\delta \mathcal{S}/\delta z=0$) to the following action:
\begin{equation}
    \mathcal{S}(z(\cdot,\cdot)) = \iint \mathcal{L}(z,z_t,z_x) \,dt\,dx
    \label{HPDEaction}
\end{equation}
where the Lagrangian density is given by
\begin{equation}
    \mathcal{L}(z,z_t,z_x)=\frac{1}{2} \langle \mathcal{M}z_t,z \rangle + \frac{1}{2} \langle \mathcal{K}z_x,z \rangle - H(z)
\end{equation}

This variational view of Hamiltonian PDEs will intrinsically guarantee the preservation of multi-symplecticity, and there will be a field generalization of Noether's theorem, which ensures conservation of momentum maps corresponding to symmetries.

Numerically, instead of discretizing the equations \eqref{HamiltonianPDE}, we prefer the approach of variational integrators because they are intrinsically multi-symplectic and therefore structure-preserving \cite{MaPaSh98, MaSh99, MaWe:01, LeMaOrWe2003}. These integrators are obtained as follows: first discretize the action \eqref{HPDEaction} using quadratures, then apply variational principle to the discrete action (which depends on finitely many arguments), and finally, solve the algebraic system obtained from the variational principle, i.e., the discrete Euler-Lagrange equations.

For an illustration, consider a nonlinear wave equation:
\begin{equation}
    u_{tt}-u_{xx}=V'(u)
    \label{nonlinearWave}
\end{equation}
with periodic boundary condition $u(x+L,t)=u(x,t)$ and compatible initial conditions $u(x,t=0)=f(x)$ and $u_t(x,t=0)=g(x)$. Suppose we are interested in the solution in a domain $[0,L]\times [0,T]$.

Rewrite the high order PDE as a system of first order PDEs (notice these covariant equations can be obtained through an intrinsic procedure, which works on manifolds as well \cite{Br06}):
\begin{eqnarray}
    v_t-w_x &=& V'(u) \\
    u_t &=& v \\
    u_x &=& w
\end{eqnarray}

The corresponding Lagrangian density is:
\begin{equation}
    \mathcal{L}=\frac{1}{2}u_t^2-\frac{1}{2}u_x^2+V(u)
\end{equation}

Using a forward time forward space approximation, we obtain the following discrete Lagrangian:
\begin{eqnarray}
    L^d_{i,j}
    &\triangleq& h_{ij}k_{ij} \left[ \frac{1}{2} \left( \frac{u_{i,j+1}-u_{i,j}}{h_{ij}} \right)^2 - \frac{1}{2} \left( \frac{u_{i+1,j}-u_{i,j}}{k_{ij}} \right)^2 + V(u_{i,j}) \right] \\
    &\approx&\int_{t_j}^{t_{j+1}=t_j+h_{ij}} dt \int_{x_i}^{x_{i+1}=x_i+k_{ij}} dx \, \left[ \frac{1}{2}u_t^2-\frac{1}{2}u_x^2+V(u) \right]
\end{eqnarray}
where space step $k_{ij}$ and time step $h_{ij}$ define a rectangular grid of size $k_{ij} \times h_{ij}$. The simplest single-scale choice would be $k_{ij}=k$ and $h_{ij}=h$ for some $k$ and $h$.

As a consequence, the continuous action $\mathcal{S}$ is approximated by a discrete action:
\begin{equation}
    \mathcal{S}_d=\sum_{\alpha=1}^{N}\sum_{\beta=1}^{M} L^d_{\alpha,\beta} \approx \mathcal{S}=\iint \mathcal{L} \,dt\,dx
\end{equation}
and Hamilton's principle of least action $\delta \mathcal{S}_d=0$ gives
\begin{equation}
    \frac{\partial}{\partial u_{i,j}} \sum_{\alpha=1}^{N}\sum_{\beta=1}^{M} L^d_{\alpha,\beta} = 0
\end{equation}
for $1\leq i \leq N$ and $1\leq j \leq M$, where $N$ and $M$ are such that $\sum_{\alpha=1}^{N}k_{\alpha\beta}=L$ for any $\beta$ and $\sum_{\beta=1}^{M}h_{\alpha\beta}=T$ for any $\alpha$.

Taking derivative with respect to $u_{i,j}$, we obtain the following discrete Euler-Lagrange equations:
\begin{equation}
    k_{ij}\frac{u_{i,j}-u_{i,j+1}}{h_{ij}}-u_{ij}\frac{u_{i,j}-u_{i+1,j}}{k_{ij}}+h_{ij}k_{ij}V'(u_{i,j})+k_{i,j-1}\frac{u_{i,j}-u_{i,j-1}}{h_{i,j-1}}-h_{i-1,j}\frac{u_{i,j}-u_{i-1,j}}{k_{i-1,j}}=0
    \label{DEL}
\end{equation}

The system of above equations is explicitly solvable when equipped with boundary conditions and initial conditions; for instance, below is a consistent discretization of the continuous version:
\begin{equation}
\begin{cases}
    u_{i,j} = u_{i+N,j},   & \forall i,j \\
    u_{i,1} = f\left(\sum_{\alpha=1}^{i} k_{\alpha 1}\right), & \forall i \\
    u_{i,2} = u_{i,1} + h_{i1} g\left(\sum_{\alpha=1}^{i} k_{\alpha 2}\right), & \forall i
\end{cases}
\end{equation}

This numerical receipt is convergent. In fact, multi-symplectic integrators obtained from variational principles can be viewed as special members of finite difference methods, whose error analysis is classical.

We wish to point out that the above procedure works for any Hamiltonian PDEs of form \eqref{HamiltonianPDE}. Also, notice that high-order derivatives are dealt with in an intrinsic way regardless of whether mesh is uniform.

\subsection{FLAVORization of multi-symplectic integrators}

Now consider a multiscale Hamiltonian PDE
\begin{equation}
    \mathcal{M}(1,\epsilon^{-1})z_t+\mathcal{K}(1,\epsilon^{-1})z_x=\nabla_z H(1,\epsilon^{-1},z)
\end{equation}

Any single-scale multi-symplectic integrator can be FLAVORized  (to achieve computational acceleration) by using  the following strategy:
(i) Use the two-scale mesh illustrated in Figure \ref{FLAVORmesh}, and (ii) Turn  off large coefficients when taking mesoscopic time-steps. Unlike FLAVORizing a general finite difference scheme, we FLAVORize the action $\mathcal{S}_d$ instead of the PDE. Specifically, choose
\begin{equation}
    \begin{cases}
    k_{ij} = K, & \forall i,j \\
    h_{ij} = h, & \forall i \text{ and odd }j \\
    h_{ij} = H-h, & \forall i \text{ and even }j
    \end{cases}
    \label{grid}
\end{equation}
and let $\epsilon^{-1}=0$ in $L^d_{i,j}$ for even $j$'s and all $i$'s, while the large value of $\epsilon^{-1}$ is kept in $L^d_{i,j}$ for odd $j$'s and all $i$'s. $h$ and $H$ correspond to a small and a mesoscopic time-step, and $K$ corresponds to a mesoscopic space-step; the same rule of thumb for choosing them in Section \ref{FDsection} applies.

After applying the discrete Hamilton's principle, the resulting discrete Euler Lagrange-equations corresponding to a multi-symplectic integrator will still be \eqref{DEL}, except that stiffness is turned off in half of the grids. Multisymplecticity is automatically gained, because the updating equations originate from a discrete variational principle \cite{MaPaSh98}.

\subsection{Example: multiscale Sine-Gordon wave equation}
Consider a specific nonlinear wave equation \eqref{nonlinearWave} in which $V(u)=-\cos(\omega u)-\cos(u)$. If $\omega=0$, this corresponds to Sine-Gordon equation, which has been studied extensively due to its soliton solutions and its relationships with quantum physics (for instance, as a nonlinear version of Klein-Gordon equation). We are interested in the case in which $\omega$ (identified with $\epsilon^{-1}$) is big, so that a separation of timescale exhibits.

Arbitrarily choose $L=2$ and use periodic boundary condition $u(x+L,t)=u(x,t)$, and let initial condition be $u(x,0)=\sin(2\pi x/L)$ and $u_t(x,0)=0$. Denote total simulation time by $T$. Use the FLAVOR mesh \eqref{grid}. In order to obtain a stable and accurate numerical solution, $k$ and $h$ have to be $o(1/\omega)$, and $K$ and $H$ need to be $o(1)$.

\begin{figure} [ht]
\begin{tabular}{cc}
\hspace{-30pt}
\includegraphics[width=0.58\textwidth]{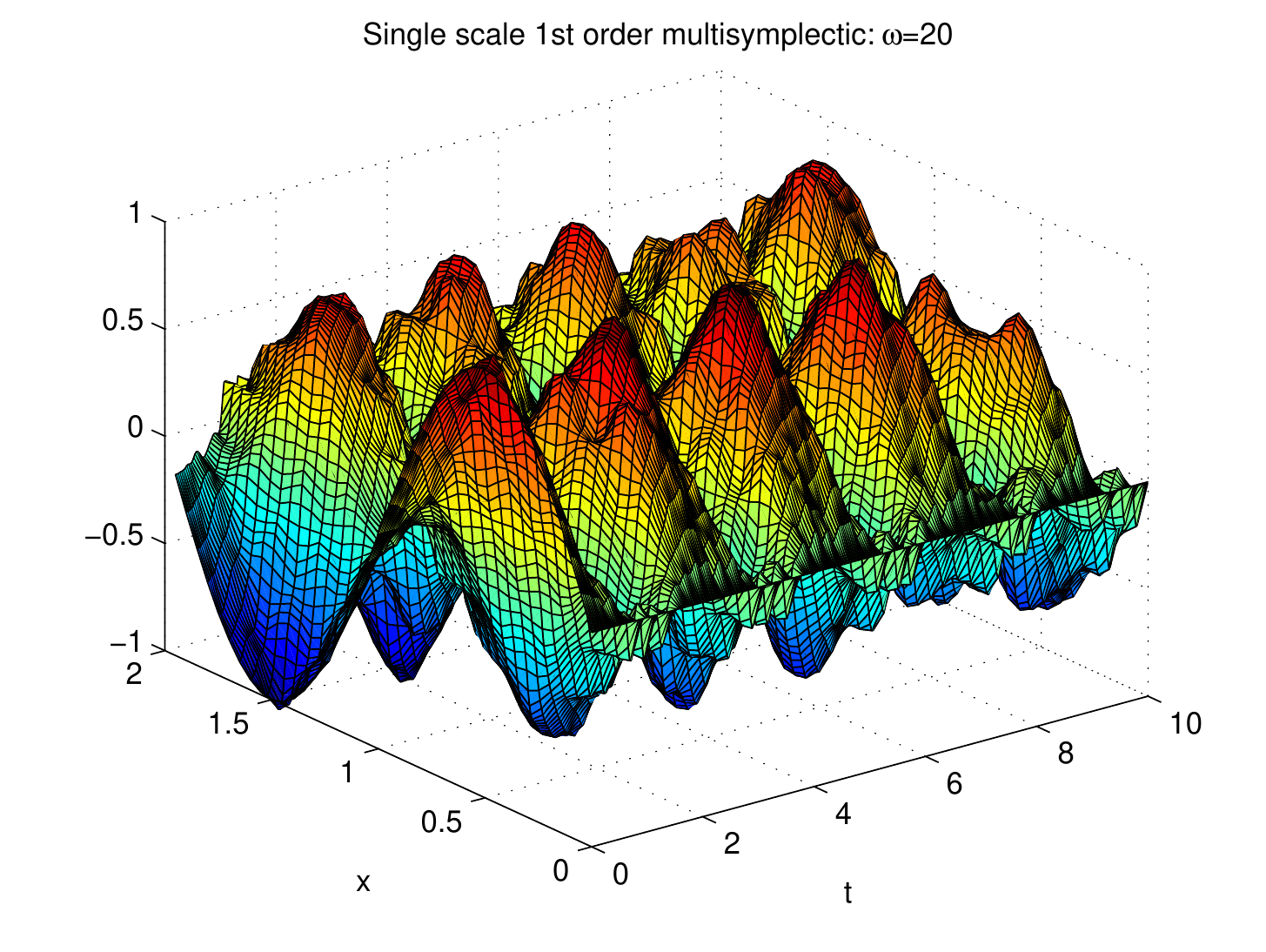}
\hspace{-20pt}
&
\hspace{-20pt}
\includegraphics[width=0.58\textwidth]{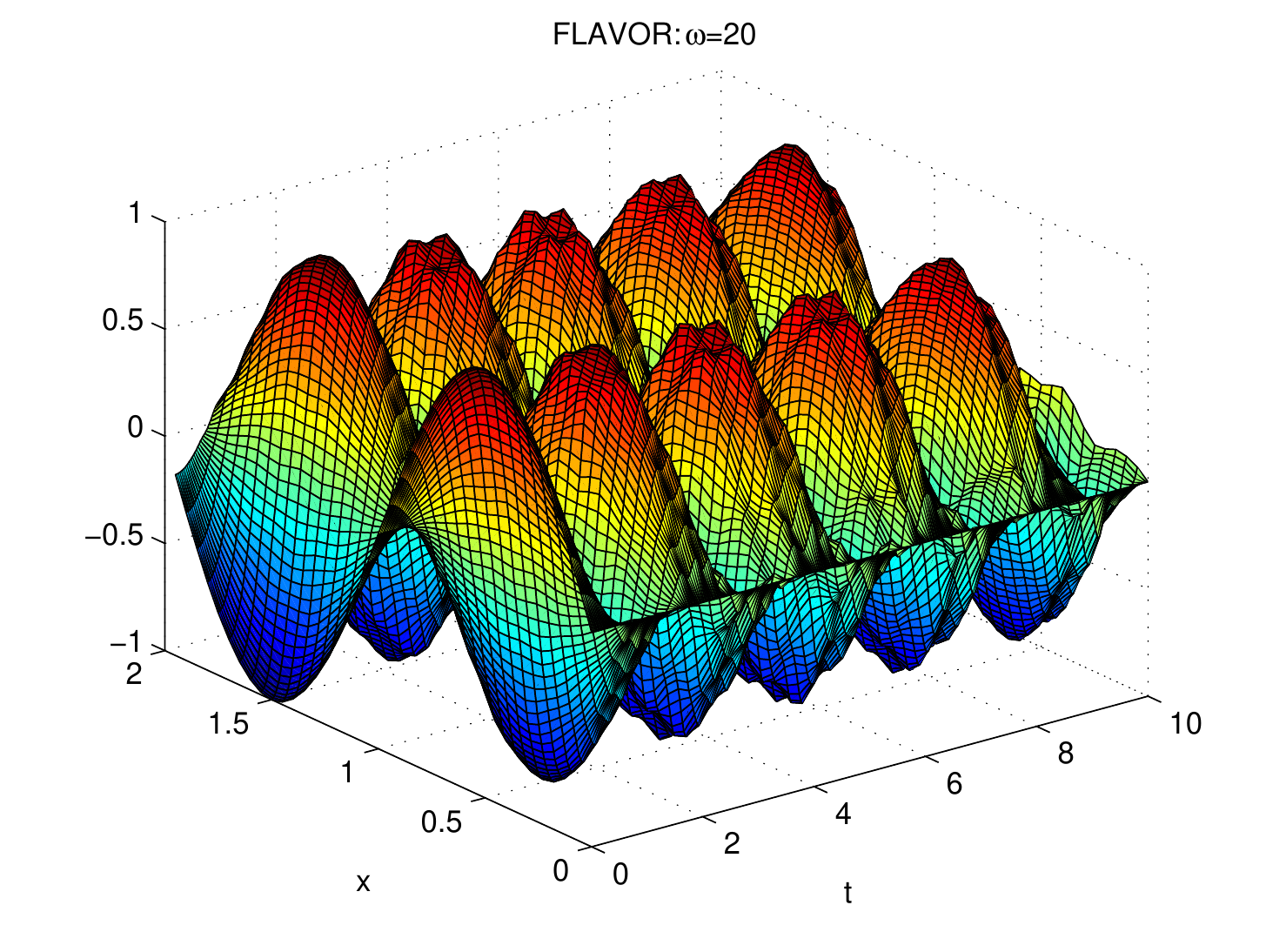}
\hspace{-30pt}
\end{tabular}
\caption{\footnotesize Numerical solutions to multiscale Sine-Gordon equation by single-scale 1st-order multi-symplectic integrator (left) and its FLAVORization (right). For clarity, the surface plots (but not simulations) use the same mesh size.}
\label{multi-symplecticResults}
\end{figure}

A comparison between the benchmark of the single-scale forward time forward space multi-symplectic integrator (Eq. \ref{DEL} with $h_{ij}=h$ and $k_{ij}=k$) and its FLAVORization (Eq. \ref{DEL} with mesh \eqref{grid} and $V'(u)=\omega\sin(\omega u)+\sin(u)$ for odd $j$ and $V'(u)=\sin(u)$ for even $j$) is presented in Figure \ref{multi-symplecticResults}. $\omega=20$, $k=L/20/\omega$ and $h=k/2$, and $K=L/40$ and $H=K/2$. It is intuitive to say that the slow process of wave propagation is well-approximated by FLAVOR, although the fast process of local fluctuation is not captured in the strong sense. Error quantification is not done, because what the slow and fast processes are is not rigorously known here. $HK/2hk=50$-fold acceleration is obtained by FLAVOR.

\begin{figure} [ht]
\center
\includegraphics[width=0.58\textwidth]{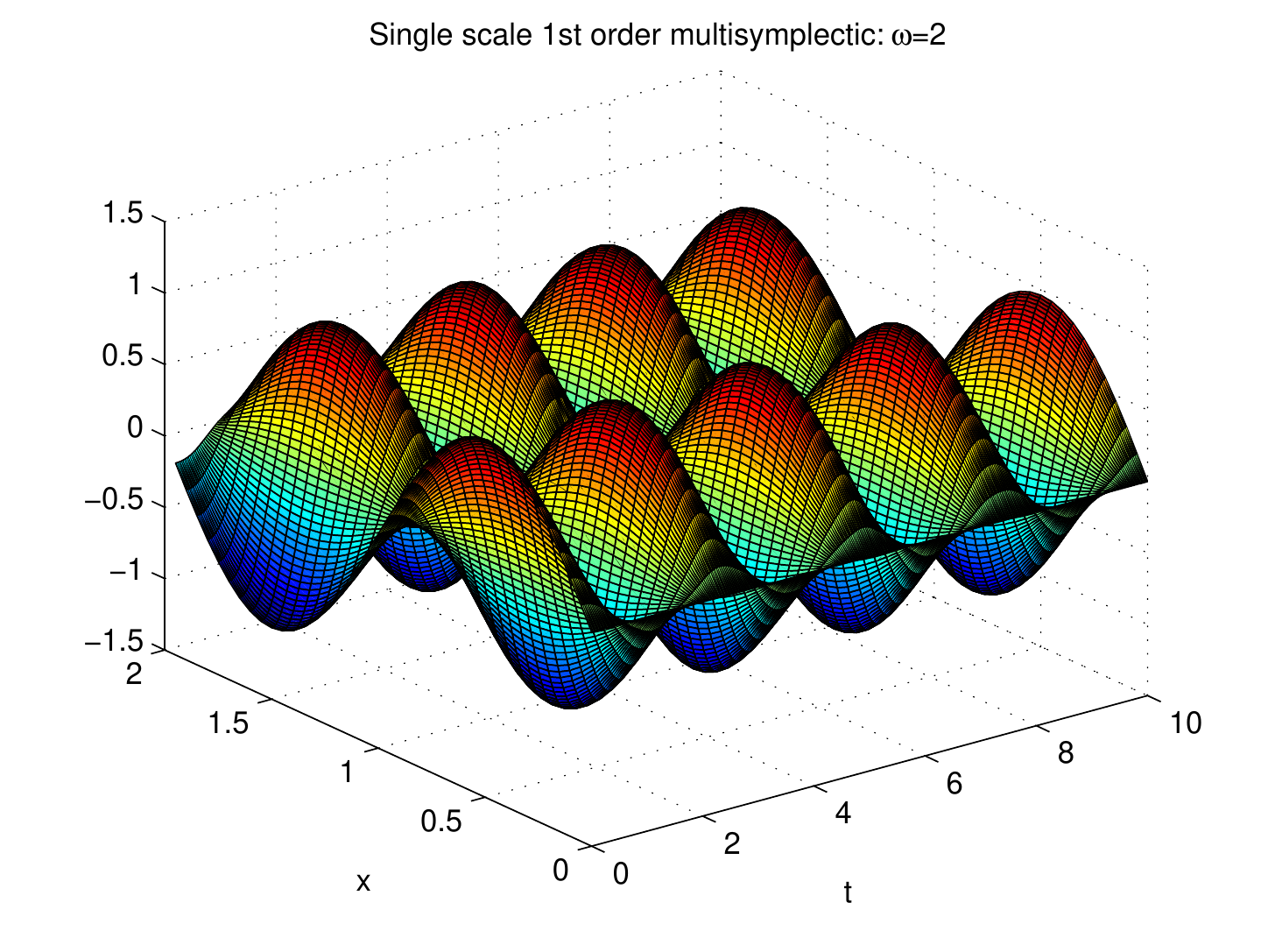}
\caption{\footnotesize Numerical solutions to multiscale Sine-Gordon equation with the `equivalent' stiffness by single-scale 1st-order multi-symplectic integrator. For clarity, the surface plot (but not the simulation) uses the same mesh size (as in Figure \ref{multi-symplecticResults}).}
\label{multi-symplecticResultsAdditional}
\end{figure}

Readers familiar with the splitting theory of ODEs \cite{McQu02} might question whether FLAVORS are equivalent to an averaged stiffness of $\tilde{\omega}=\omega \frac{h}{H}$ (which corresponds $\tilde{\omega}=2$ in the numerical experiment described above). The answer is no, because the equivalency given by the splitting theory is only local. In fact, the same single-scale forward time forward space multi-symplectic integration of the case $\omega=2$ is shown in Figure \ref{multi-symplecticResultsAdditional}, which is clearly distinct from the FLAVOR result in Figure \ref{multi-symplecticResults}. Moreover, because of the $e^{C \omega T}$ error term, changing stiffness alone will not result in a converging method (and result in a $\mathcal{O}(1)$ error on slow variables).

\section{Pseudospectral methods}\label{secpseudospectral}
\subsection{Single-scale method}
Consider a PDE
\begin{equation}
    u_t(x,t)=\mathcal{L} u(x,t)
\end{equation}
with periodic boundary condition $u(x,t)=u(x+L,t)$ and initial condition $u(x,0)=f(x)$, where $\mathcal{L}$ is a differential operator involving only spatial derivatives.

The Fourier collocation method approximates the solutions by the truncated Fourier series:
\begin{equation}
    u_N(x,t)=\sum_{|n|\leq N/2} a_n(t) e^{in2\pi x/L}
    \label{Fourier}
\end{equation}
and solves for $a_n(t)$'s by requiring the PDE to hold at collocation points $y_j$:
\begin{equation}
    \partial_t u_N(y_j,t)-\mathcal{L} u_N(y_j,t)=0
    \label{pseudoODEs}
\end{equation}
This yields $N$ ODEs, which can be integrated by any favorite ODE solver. Of course, specific choices of collocations points will affect the numerical approximation. Oftentimes, the simplest choice of $y_j=Lj/N,j=0,\ldots,N-1$ is used, and in this case, the method is also called a pseudospectral method.
We refer to \cite{HeGoGo07}  for additional details on Fourier collocation methods.  It is worth mentioning that pseudospectral methods can also be multi-symplectic when applied to Hamiltonian PDEs \cite{Ch06}.

\subsection{FLAVORization of pseudospectral methods}
When the PDE is stiff (for instance, when $\mathcal{L}$ contains a large parameter $\epsilon^{-1}$), FLAVORS can be employed to integrate the  stiff ODEs (which will still contain $\epsilon^{-1}$) resulting from a pseudospectral discretization.

Similarly, for the FLAVORization of a pseudospectral method, it is sufficient to choose $N \gg L$ instead of $N \gg \epsilon^{-1} L$, i.e., the space-step can be coarse ($K=o(1)$). For time stepping, alternatively switching between $h=o(\epsilon)$ and $H-h$ for a mesoscopic $H=o(1)$ is again needed, and stiffness has to be turned off over the mesoscopic step of $H-h$. In a sense, we are still using the same FLAVOR `mesh' (Figure \ref{FLAVORmesh}), except that here we do not discretize space, but instead truncate Fourier series to resolve the same spatial grid size.

\subsection{Example: a slow process driven by a non-Dirac fast process}
Consider the following system of PDEs
\begin{equation}
    \begin{cases}
        u_t+u_x-q^2=0 \\
        q_t+q_x-p=0 \\
        p_t+p_x+\omega^2 q=0
    \end{cases}
    \label{madeupExample}
\end{equation}
with periodic boundary conditions $u(x,t)=u(x+L,t)$, $q(x,t)=q(x+L,t)$, and $p(x,t)=p(x+L,t)$, and initial conditions $u(x,0)=f^u(x)$, $q(x,0)=f^q(x)$, and $p(x,0)=f^p(x)$. The integration domain is restricted to $[0,T]\times[0,L]$. The stiffness $\epsilon^{-1}$ is identified with $\omega^2$.
We choose the initial condition of $f^u(x)=f^q(x)=\cos(2\pi x/L)$ and $f^p(x)=0$.

In this system, $q$ and $p$ correspond to a fast process, which is a field theory version of a harmonic oscillator with high frequency $\omega$. $u$ is a slow process, into which energy is pumped by the fast process in a non-trivial way.

We have chosen to FLAVORize \eqref{madeupExample} because it does not fall into the (simpler) category of systems with fast processes converging towards Dirac (single point support) invariant distributions \cite{FiJi10}.

We use the classical 4th order Runga-Kutta scheme (see, for instance, \cite{MR1227985}) for the (single-step) time integration of the pseudospectrally discretized system of ODEs \eqref{pseudoODEs}. Write $\phi_h^{\omega^2}: \tilde{a}_n^{u,q,p}(t) \mapsto \tilde{a}_n^{u,p,q}(t+h)$ its
numerical flow  over a microscopic time step $h$ (consisting of four sub-steps), where $\tilde{a}_n^{u,q,p}(t)$ are numerical approximations to the Fourier coefficients in \eqref{Fourier},  for the unknowns $u$, $q$ and $p$ at an arbitrary time $t$. Then, the corresponding FLAVOR update over a mesoscopic time step $H$ will be $\phi_{H-h}^{0} \circ \phi_h^{\omega^2}$, which consists of eight sub-steps.

\begin{figure} [ht]
\begin{tabular}{cc}
\hspace{-30pt}
\includegraphics[width=0.58\textwidth]{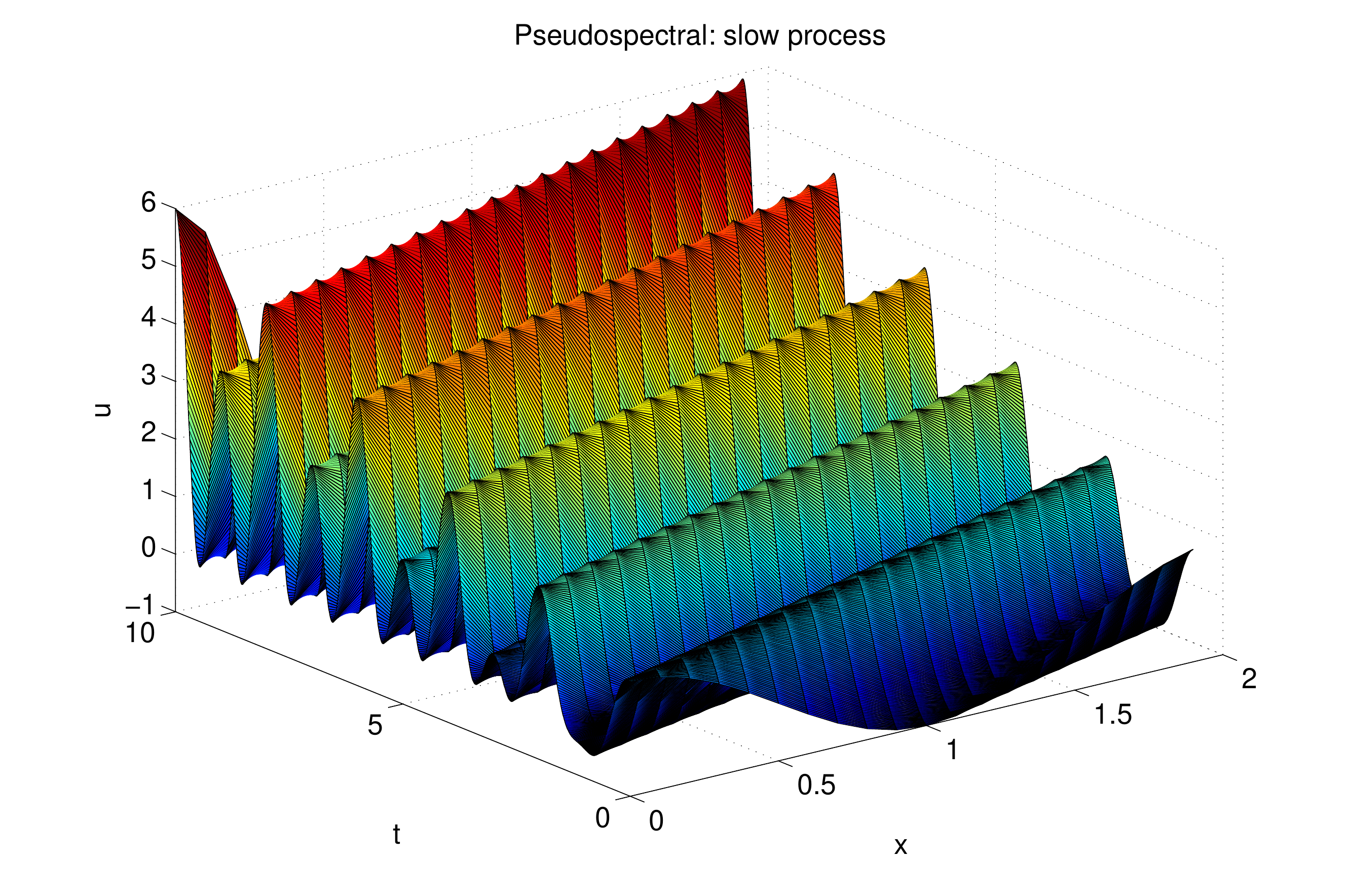}
\hspace{-20pt}
&
\hspace{-20pt}
\includegraphics[width=0.58\textwidth]{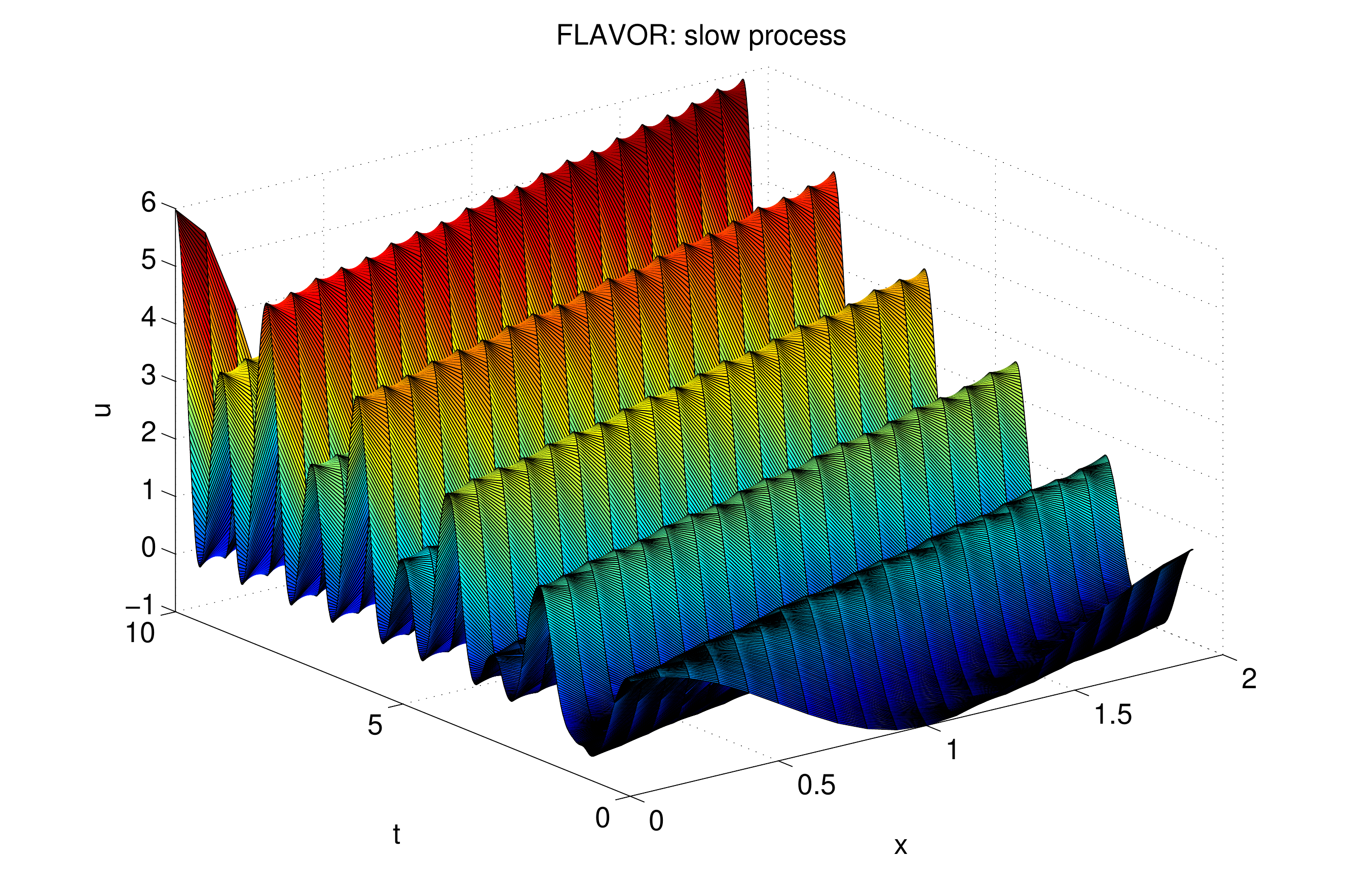}
\hspace{-30pt}
\end{tabular}
\caption{\footnotesize Single-scale (left) and multiscale pseudospectral (right) integrations of slow $u$ in system \eqref{madeupExample}. Plotting mesh for the single-scale simulation is coarser than its computation mesh.}
\label{pseudospectralResultsSlow}
\end{figure}

\begin{figure} [ht]
\begin{tabular}{cc}
\hspace{-30pt}
\includegraphics[width=0.58\textwidth]{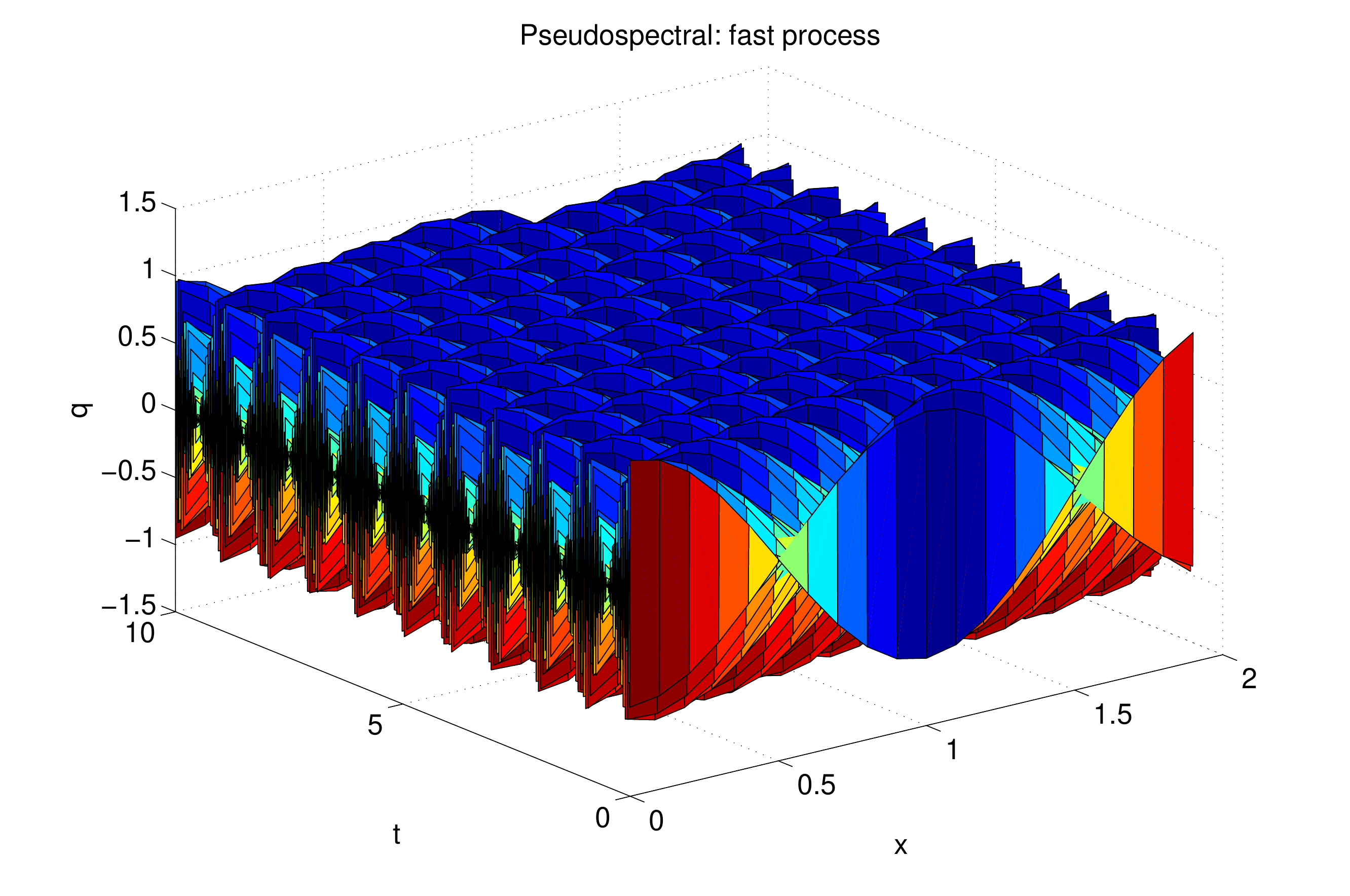}
\hspace{-20pt}
&
\hspace{-20pt}
\includegraphics[width=0.58\textwidth]{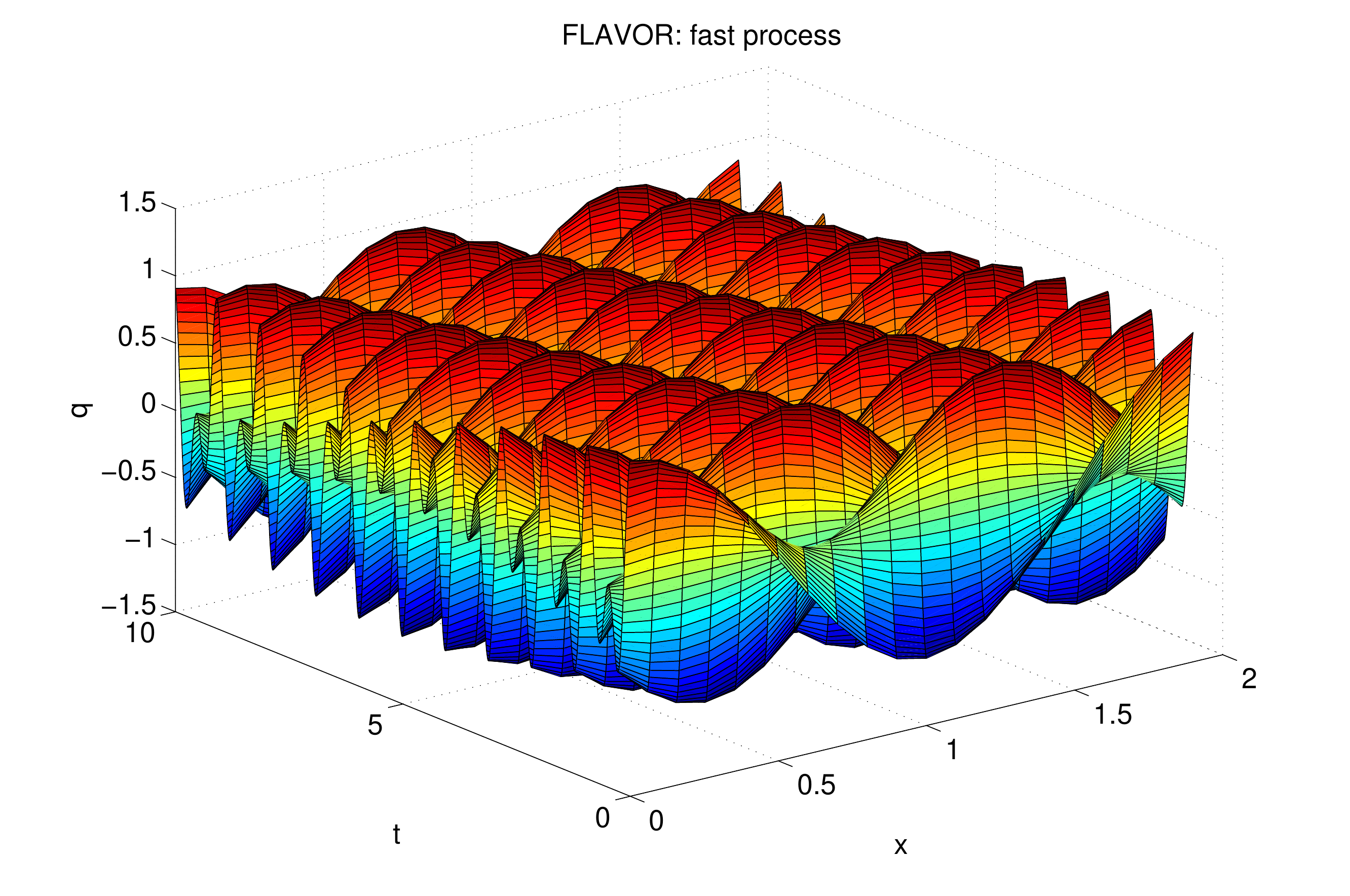}
\hspace{-30pt}
\end{tabular}
\caption{\footnotesize Single-scale (left) and multiscale pseudospectral (right) integrations of fast $q$ in system \eqref{madeupExample}. Plotting mesh for the single-scale simulation is coarser than its computation mesh. The same color does not indicate the same value in these two plots.}
\label{pseudospectralResultsFast}
\end{figure}

We present in Figure \ref{pseudospectralResultsSlow} and Figure \ref{pseudospectralResultsFast} a comparison between the benchmark of single-scale pseudospectral simulation and its FLAVORization. It can be seen that the slow process of $u$ is captured in strong (point-wise) sense, whereas the fast process of $q$ is only approximated in a weak sense (i.e. as a measure, in the case, wave shape and amplitude are correct, but not the period). We choose $L=2$, $T=10$ and $\omega=1000$. The  single-step integration uses $N=20$ and $h=0.1/\omega$ (notice that this is already beyond the stability/accuracy region of a single-scale finite difference, since the space step does not depend on $1/\omega$; the spectral method is more stable/accurate for a large space-step), and FLAVOR uses $N=20$, $h=1/\omega^2$ and $H=0.01$. $H/2h=50$-fold acceleration is achieved by FLAVOR.

\section{Convergence analysis}
\label{SectionErrorAnalysis}
\subsection{Semi-discrete system}\label{jhsjgfdhgdfdhgf3}
All  FLow AVeraging integratORS  described in previous sections are illustrations of the following (semi-discrete) strategy:
 first, space is discretized or interpolated; next, spatial differential operators are approximated by algebraic functions of finitely many spatial variables; finally, the resulting system of ODEs  is  numerically integrated by a corresponding ODE-FLAVOR \cite{FLAVOR10}. In this section, we will use the semi-discrete ODE system as an intermediate link to demonstrate that these PDE-FLAVORS are convergent to the exact PDE solution  under reasonable assumptions (in a strong sense with respect to (possibly hidden) slow variables and in the sense of measures with respect to fast variables).

More precisely, consider a spatial mesh (vector) $\mathcal{M}^S=[x_1,x_2,\ldots]$, a temporal mesh (vector) $\mathcal{M}^T=[t_1,t_2,\ldots]$, and a domain mesh (matrix) $\mathcal{M}=\mathcal{M}^S \times \mathcal{M}^T$. Examples of these meshes include the FLAVOR mesh $\mathcal{M}^S=[K,2K,\ldots,NK]$ and $\mathcal{M}^T=[h,H,H+h,2H,\dots,(M-1)H,(M-1)H+h,MH]$, and a usual single-scale (step) integration mesh $\mathcal{M}^S=[k,2k,\ldots,L]$ and $\mathcal{M}^T=[h,2h,\dots,T]$ (recall the domain size is $L=NK$ by $T=MH$). We will use the FLAVOR mesh throughout this section. We will compare the solution of the PDE \eqref{1stOrderPDE} with the solution obtained with the FLAVOR strategy at these discrete points.

For simplicity, assume the PDE of interest is 1st-order in time derivative:
\begin{equation}
    u_t(x,t)=F(1,\epsilon^{-1},x,t,u(x,t),u_x(x,t),\ldots)
    \label{1stOrderPDE}
\end{equation}
Observe that a PDE \eqref{generalPDE} with higher-order time derivatives can be written as a system of 1st-order (in time derivatives) PDEs.

Now  consider a consistent discretization of PDE \eqref{1stOrderPDE} with space step $K$ and time step $h$ (we refer to   Page 20 of \cite{St04} for a definition of the notion of consistency, which intuitively means vanishing local truncation error). Letting $h\downarrow 0$ in this discretization, we obtain a semi-discrete system (continuous in time and discrete in space). This semi-discrete system is denoted by the following system of ODEs, with approximated spatial derivatives:
\begin{equation}
\begin{cases}
    \dot{u}_1(t)=f_1(u_1,u_2,\ldots,u_N,\epsilon^{-1},t) \\
    \dot{u}_2(t)=f_2(u_1,u_2,\ldots,u_N,\epsilon^{-1},t) \\
    \quad \cdots    \\
    \dot{u}_N(t)=f_N(u_1,u_2,\ldots,u_N,\epsilon^{-1},t)
\end{cases}
\label{semi-discrete}
\end{equation}
Assuming  existence and uniqueness of an exact $\mathcal{C}^1$ strong solution $u$ to the PDE \eqref{1stOrderPDE}, and writing $u(\mathcal{M}^S_i,t)$ its values at the spatial discretization points, we define for each $i$ the following remainder:
\begin{equation}\label{eqreminder}
    \mathcal{R}_i(\epsilon^{-1},t) \triangleq \frac{\partial u}{\partial t}(\mathcal{M}^S_i,t) - f_i(u(\mathcal{M}^S_1,t),u(\mathcal{M}^S_2,t),\ldots,u(\mathcal{M}^S_N,t),\epsilon^{-1},t)
\end{equation}
which is a real function of $t$ indexed by $\epsilon^{-1}$.

Then, $u_i(t)$ approximates the exact solution $u(\mathcal{M}^S_i,t)$ evaluated at grid points in the sense that these remainders vanish as $\epsilon^{-1} K\downarrow 0$ (where $K:=\mathcal{M}^S_i-\mathcal{M}^S_{i-1}$):
\begin{Lemma}
    Assume that $F$ in \eqref{1stOrderPDE} satisfies
    \begin{equation}
        |F(1,\epsilon^{-1},x,t,u(x,t),u_x(x,t),\ldots)| \leq (1+\epsilon^{-1}) |F(1,1,x,t,u(x,t),u_x(x,t),\ldots)|
        \label{g28094yt2puboiadjglbhga}
    \end{equation}
    Assume that the $f_i$ in \eqref{semi-discrete} satisfy similar inequalities.
    Then, there exists a constant $C_i$ independent from $\epsilon$, $h$, $H$ or $K$, such that for bounded $t$ and $u$
    \begin{equation}
        |\mathcal{R}_i(\epsilon^{-1},t)|\leq (1+\epsilon^{-1}) C_i K
        \label{spatialError}
    \end{equation}
    \label{LemmaSpatialError}
\end{Lemma}

\begin{Remark}
    \eqref{g28094yt2puboiadjglbhga} is true, for instance, in cases where
    \begin{equation}
        F(1,\epsilon^{-1},x,t,u(x,t),\ldots)=F_0(x,t,u(x,t),\ldots)+\epsilon^{-1}F_1(x,t,u(x,t),\ldots).
    \end{equation}
\end{Remark}

\begin{proof}
    The linear scaling with $K$ in \eqref{spatialError} immediately follows from the definition of consistency, and the parameter $1+\epsilon^{-1}$ in \eqref{spatialError} has its origin \eqref{g28094yt2puboiadjglbhga}.
\end{proof}

\begin{Remark}
    The consistency of finite difference methods can be easily shown using Taylor expansions. For instance, applying a Taylor expansion to the solution of  $u_t-\epsilon^{-1} u_x=a(u)$ leads to
    \begin{equation}
    \begin{split}
        u(iK,(j+1)h)=&u(iK,jh)+h\Big(\epsilon^{-1} \big(\frac{u((i+1)K,jh)-u(iK,jh)}{K}\\&+\mathcal{O}(K)\big)+a(u(iK,jh))\Big)+\mathcal{O}(h^2)
        \end{split}
    \end{equation}
  which implies
    \begin{equation}
        \frac{\partial}{\partial t} u(iK,t)=\epsilon^{-1} \frac{u((i+1)K,t)-u(iK,t)}{K}+a(u(iK,t))+\epsilon^{-1}\mathcal{O}(K)
    \end{equation}
   and naturally establishes  the correspondence of  $f_i(u_1,\ldots,u_N,\epsilon^{-1},t)=\epsilon^{-1} \frac{u_{i+1}(t)-u_i(t)}{K}+a(u_i(t))$ and $\mathcal{R}_i=\epsilon^{-1}\mathcal{O}(K)$  for a 1st-order finite difference scheme.
    Notice that the remainders are still stiff, but we will see later that this is not a problem, since they can be handled by ODE-FLAVORs.
 The consistency of pseudospectral method can be shown similarly using Fourier analysis.
\end{Remark}

With $\mathcal{R}_i$  defined in \eqref{eqreminder}, consider the following system of ODEs:
\begin{equation}
\begin{cases}
    \dot{u}_1(t)=f_1(u_1,u_2,\ldots,u_N,\epsilon^{-1},t)+\mathcal{R}_1(\epsilon^{-1},t) \\
    \quad \cdots    \\
    \dot{u}_N(t)=f_N(u_1,u_2,\ldots,u_N,\epsilon^{-1},t)+\mathcal{R}_N(\epsilon^{-1},t)
\end{cases}
\label{semi-discrete_exact}
\end{equation}
with initial condition $u_i(0)=u(\mathcal{M}^S_i,0)$. Obviously, its solution $(u_i(t))_{1\leq i \leq N}$ is the exact PDE solution sampled at spatial grid points, i.e., $u_i(t)=u(\mathcal{M}^S_i,t)$.

We will now establish the accuracy of  PDE-FLAVOR by showing that an ODE-FLAVOR integration of \eqref{semi-discrete_exact} leads to an accurate approximation of $(u_i(t))_{1\leq i \leq N}$. Since  space (with fixed width $L$) is discretized  by $N$ grid points, we use the following (normalized by $N$) norm in our following discussion (suppose $v_i(t)=v(\mathcal{M}^S_i,t)$ for a function $v$):
\begin{equation}
    \left\| [v_1(t),v_2(t),\ldots,v_N(t)] \right\| \triangleq \frac{1}{N}\left\| [v_1(t),v_2(t),\ldots,v_N(t)] \right\| _1
    \label{DefinitionNorm}
\end{equation}
Observe that if $v(\cdot,t)$ is Riemann integrable, then
\begin{equation}
    \lim_{K\downarrow 0} \left\| [v(\mathcal{M}^S_1, t),v(\mathcal{M}^S_2,t),\ldots,v(\mathcal{M}^S_N, t)]
    \right\| \rightarrow \frac{1}{L}\left\| v(\cdot,t) \right\|_{\mathcal{L}^1} \qquad (\text{recall }L=NK \text{ is fixed}),
\end{equation}
and hence the norm \eqref{DefinitionNorm} does not blow up or vanish as $N\rightarrow \infty$.

\subsection{Sufficient conditions for convergence, ODE-FLAVORS, and two-scale convergence of PDE-FLAVORS}
\label{SectionLala}

We will now prove the accuracy of PDE-FLAVORs under the assumption of existence of (possibly hidden) slow and locally ergodic fast variables.
The convergence of PDE-FLAVORs will be expressed using the notion of two-scale flow convergence introduced in \cite{FLAVOR10} (corresponding to a strong convergence with respect to slow variables and weak one with respect to fast ones).

\begin{Condition}
    Assume that the ODE system \eqref{semi-discrete_exact} satisfies the following conditions:
    \begin{enumerate}
    \item
        (Existence of hidden slow and fast variables): There exists a (possibly time-dependent) diffeomorphism $\eta^t:[u_1(t),\ldots,u_N(t)] \mapsto [x(t),y(t)]$ from $\mathbb{R}^N$ onto $\mathbb{R}^{N-p}\times \mathbb{R}^p$ with uniformly bounded $C^1, C^2$ derivatives with respect to $u_i$'s and $t$, and such that for all $\epsilon>0$, $(x(t),y(t))$ satisfies
        \begin{equation}
        \begin{cases}
            \dot{x}(t)&=f(x(t),y(t),t) \\
            \dot{y}(t)&=\epsilon^{-1} g(x(t),y(t),t)
        \end{cases} \quad ,
        \end{equation}
        where $f$ and $g$ have bounded $C^1$ derivatives with respect to $x$, $y$ and $t$.
    \item
        (Local ergodicity of vast variables):   There exists a family of probability measures $\mu^t(x,dy)$ on $\mathbb{R}^p$ indexed by $x\in\mathbb{R}^{N-p}$ and $t\in\mathbb{R}$, and a family of positive functions $T \mapsto E^t(T)$ satisfying $\lim_{T\rightarrow \infty}E^t(T)=0$ for all bounded $t$, such that for all $x_0,y_0,t_0,T$ bounded and $\phi$ uniformly bounded and Lipschitz, the solution to
        \begin{equation}
            \dot{Y}_t=g(x_0,Y_t,t_0)\quad \quad Y_0=y_0
        \end{equation}
        satisfies
        \begin{equation}
            \Big|\frac{1}{T}\int_0^T \phi(Y_s)ds-\int_{\mathbb{R}^p} \phi(y)\mu^{t_0}(x_0,dy)\Big|\leq \chi^{t_0} \big(\|(x_0,y_0)\|\big) E^{t_0}(T) (\|\phi\|_{L^\infty}+\|\nabla \phi\|_{L^\infty} )
        \end{equation}
        where $r \mapsto \chi^{t_0}(r)$ is bounded on compact sets, and $\mu^t$ has bounded derivative with respect to $t$ in total variation norm.
    \end{enumerate}
    \label{ContinuousFLAVORcondition}
\end{Condition}

Under Conditions \ref{ContinuousFLAVORcondition}, the computation of the solution of PDE \eqref{1stOrderPDE} can be accelerated by applying the FLAVOR strategy to  a single-scale time integration of  its semi-discretization \eqref{semi-discrete}.

Write $\Phi_{t,t+\tau}^{\alpha}$ the numerical flow of a given (legacy) ODE integrator for \eqref{semi-discrete}:
\begin{equation}
    \Phi_{t,t+\tau}^{\alpha}: [\tilde{u}_1(t),\ldots, \tilde{u}_N(t)] \mapsto [\tilde{u}_1(t+\tau),\ldots, \tilde{u}_N(t+\tau)] \, ,
\label{eqlegacyint}
\end{equation}
where $\tilde{u}_i(s)$ approximates $u_i(s)$ for all $s$, $\tau$ is the integration time step, and $\alpha$ is a controllable parameter that replaces the stiff parameter $\epsilon^{-1}$ in \eqref{semi-discrete} and takes values of $\epsilon^{-1}$ (stiffness `on') or $0$ (stiffness `off').

\begin{Definition}[ODE-FLAVORS]
    The FLow AVeraging integratOR associated with $\Phi$ is defined as the algorithm simulating the process:
    \begin{align}\label{FLAVOR_ODE}
        & \quad [\bar{u}_1(t),\ldots,\bar{u}_N(t)] \nonumber\\
        &= \big( \Phi^0_{(k-1)H+h,kH}\circ \Phi^{\frac{1}{\epsilon}}_{(k-1)H,(k-1)H+h} \big) \circ \cdots \circ \big( \Phi^0_{H+h,2H}\circ \Phi^{\frac{1}{\epsilon}}_{H,H+h} \big) \circ \big( \Phi^0_{h,H}\circ \Phi^{\frac{1}{\epsilon}}_{0,h} \big) ([u_1(0),\ldots,u_N(0)])
    \end{align}
    where (the number of steps) $k$ is a piece-wise constant function of $t$ satisfying $kH \leq t <(k+1)H$,
     $h$ is a microscopic time step resolving the fast timescale ($h \ll\epsilon$), $H$ is a mesoscopic time step independent of the fast timescale satisfying $h \ll\epsilon \ll H \ll 1$ and
    \begin{equation}\label{eqlimits}
        (\frac{h}{\epsilon})^2\ll H \ll \frac{h}{\epsilon}
    \end{equation}
    \label{FLAVOR_ODE_integrator}
\end{Definition}

\begin{Condition}\label{ConditionTemporalConsistency}
    Consider the legacy ODE integrator with one-step update map $\Phi_{t,t+\tau}^{\alpha}$ introduced in \eqref{eqlegacyint}. Suppose there exists constants $C>0$ and $H_0>0$ independent of $N$ and $\alpha$, such that for any $\tau \leq H_0 \min(1/\alpha, 1)$ and bounded vector $[u_1,\ldots,u_N]$,
    \begin{align}
        & \| \Phi_{t,t+\tau}^{\alpha}(u_1,\ldots,u_N)-[u_1,\ldots,u_N]-\tau[f_1(u_1,\ldots,u_N,\alpha,t),\ldots \nonumber \\
        & \qquad \ldots, f_N(u_1,\ldots,u_N,\alpha,t)] \| \leq C \tau^2 (1+\alpha)^2 \, ,
        \label{TemporalConsistency}
    \end{align}
\end{Condition}
Condition \ref{ConditionTemporalConsistency} corresponds to the assumption that  the integrator $\Phi_{t,t+\tau}^{\alpha}$ is consistent for \eqref{semi-discrete}.

Observe that we are integrating \eqref{semi-discrete} but not \eqref{semi-discrete_exact}, since the remainders $\mathcal{R}_i$'s are a-priori unknown unless the exact PDE solution is known. However, the following lemma implies that the FLAVORization of this integration is in fact convergent to the solution of \eqref{semi-discrete_exact}, even though $\mathcal{R}_i$'s are possibly stiff.

\begin{Lemma}
    Assume that $\Phi_{t,t+\tau}^{\alpha}$, introduced in \eqref{eqlegacyint},  satisfies Condition \ref{ConditionTemporalConsistency}. Let $h$ and $H$ be the time steps used in the FLAVORization \ref{FLAVOR_ODE_integrator}.
     If $h \ll \epsilon$, $H \ll h/\epsilon$, and $K= \mathcal{O}(H)$, then
    \begin{align}
        & \| \Phi_{t,t+\tau}^{\alpha}(u_1,\ldots,u_N)-[u_1,\ldots,u_N]-\tau[f_1(u_1,\ldots,u_N,\alpha,t)+\mathcal{R}_1(\alpha,t),\ldots \nonumber \\
        & \qquad \ldots, f_N(u_1,\ldots,u_N,\alpha,t)+\mathcal{R}_N(\alpha,t)] \| \leq C \tau^2 (1+\alpha)^2
        \label{TemporalConsistencyExact}
    \end{align}
    where $\tau=h$ when $\alpha=\epsilon^{-1}$ and $\tau=H-h$ when $\alpha=0$.
    \label{LemmaConsistentTemporal}
\end{Lemma}
\begin{proof}
    By Condition \ref{ConditionTemporalConsistency}, we have
    \begin{align}
        & \| \Phi_{t,t+\tau}^{\alpha}(u_1,\ldots,u_N)-[u_1,\ldots,u_N]-\tau[f_1(u_1,\ldots,u_N,\alpha,t),\ldots \nonumber \\
        & \qquad \ldots, f_N(u_1,\ldots,u_N,\alpha,t)] \| \leq C \tau^2 (1+\alpha)^2
    \end{align}
    for any $\tau\leq \min(1/\alpha,1)H_0$.
    In addition, Lemma \ref{LemmaSpatialError} gives a bound on the remainders: when $\alpha=\epsilon^{-1}$, there exists a constant $\tilde{C}>0$ independent of $N$ and $\epsilon^{-1}$, such that for all $i$,
    \begin{equation}
        |\tau \mathcal{R}_i(\epsilon^{-1},t)|\leq \tau \tilde{C} K \epsilon^{-1}
    \end{equation}
    Because we use $\tau=h$ in this case and $K\ll \epsilon^{-1} \tau$, the above is bounded by $\tau \tilde{C} (\hat{C} \epsilon^{-1}\tau) \epsilon^{-1} \leq C \tau^2(1+\alpha)^2$ for some constants $\hat{C} \ll 1$ and $C=\tilde{C}\hat{C}$.
    When $\alpha=0$ on the other hand, there exists a constant $\tilde{C}>0$ such that for all $i$
    \begin{equation}
        |\tau \mathcal{R}_i(\epsilon^{-1},t)|\leq \tau \tilde{C} K
    \end{equation}
    Because $K=\mathcal{O}(H)$ and we use $\tau=H-h=\mathcal{O}(H)$ in this case, the above is bounded by $\tau \tilde{C} \hat{C} \tau \leq C \tau^2(1+\alpha)^2$ for some constants $\hat{C}$ and we let $C=\tilde{C}\hat{C}$.
    Notice that the value of $K$ is fixed in both cases but $\tau$ has different values: the flow map used in FLAVOR associated with $\alpha=0$ is the one with mesoscopic step $\Phi_{t+h,t+H}^0$, i.e., $\tau=H-h$; when $\alpha=\epsilon^{-1}$ on the other hand, the flow map is $\Phi_{t,t+h}^{\epsilon^{-1}}$ and $\tau=h$.
    Finally, the triangle inequality gives
    \begin{align}
        & \| \Phi_{t,t+\tau}^{\alpha}(u_1,\ldots,u_N)-[u_1,\ldots,u_N]-\tau[f_1(u_1,\ldots,u_N,\alpha,t)+\mathcal{R}_1(\alpha,t),\ldots \nonumber \\
        & \ldots, f_N(u_1,\ldots,u_N,\alpha,t)+\mathcal{R}_N(\alpha,t)] \| \leq \| \Phi_{t,t+\tau}^{\alpha}(u_1,\ldots,u_N)-[u_1,\ldots,u_N]- \nonumber \\
        & \tau[f_1(u_1,\ldots,u_N,\alpha,t), \ldots, f_N(u_1,\ldots,u_N,\alpha,t)] \| + \frac{1}{N}\sum_{i=1}^N |\tau \mathcal{R}_i(\alpha,t)| \leq 2 C \tau^2 (1+\alpha)^2
         \, , \label{ofgqr7yhphgrqeui}
    \end{align}
    which finished the proof after absorbing the coefficient  $2$  into $C$.
\end{proof}

We also need the usual regularity and stability assumptions to prove the accuracy of FLAVORS for \eqref{semi-discrete_exact}.

\begin{Condition}Assume that
    \begin{enumerate}
    \item
        $f_1, f_2, \ldots, f_N$ are Lipschitz continuous.
    \item
        For all bounded initial condition $[u_1(0),\ldots,u_N(0)]$'s, the exact trajectories \\ $([u_1(t),\ldots,u_N(t)])_{0\leq t\leq T}$ (i.e., solution to \eqref{semi-discrete_exact}) are uniformly bounded in $\epsilon$.
    \item
        For all bounded initial condition $[u_1(0),\ldots,u_N(0)]$'s, the numerical trajectories  $([\bar{u}_1(t),\ldots,\bar{u}_N(t)])_{0\leq t\leq T}$ (defined by \eqref{FLAVOR_ODE}) are uniformly bounded in $\epsilon$, $0< H \leq H_0$, $h \leq \min(H_0 \epsilon, H)$.
    \end{enumerate}
    \label{RegularityCondition}
\end{Condition}

The following theorem shows the two-scale flow convergence (strong on slow variables $x$ and in the sense of measures on fast ones $y$, see \cite{FLAVOR10}) of FLAVORs under the above conditions.

\begin{Theorem}\label{ODEFLAVORTheorem}
    Consider FLAVOR trajectories in Definition \ref{FLAVOR_ODE_integrator}. Under Conditions \ref{ContinuousFLAVORcondition}, \ref{ConditionTemporalConsistency} and \ref{RegularityCondition}, there exist $C>0$, $\hat{C}>0$ and $H_0>0$ independent from $\epsilon^{-1}$ and $N$, such that for $K/\hat{C} < H <H_0$, $h <H_0 \epsilon $ and $t>0$,
    \begin{equation}
        \|x(t)-[\eta^t]^x(\bar{u}_1(t),\ldots,\bar{u}_N(t))\|\leq C e^{C t} \chi_1(u_1(0),\ldots,u_N(0),\epsilon,H,h)
        \label{jkdsgkdjshgdsdkjghe}
    \end{equation}
    and for all bounded and uniformly Lipschitz continuous test functions $\varphi:\mathbb{R}^N\mapsto \mathbb{R}$,
    \begin{align}
        & \left| \frac{1}{\Delta t}\int_{t}^{t+\Delta t}\varphi([\bar{u}_1(s),\ldots,\bar{u}_N(s)])\,ds  -\int_{\mathbb{R}^p} \varphi([\eta^t]^{-1}(x(t) ,y))\mu^t(x(t) ,dy) \right| \nonumber \\
        & \qquad \qquad \qquad  \leq \chi_2(u_1(0),\ldots,u_N(0),\epsilon,H,h,\Delta t,t) (\|\varphi\|_{L^\infty}+\|\nabla \varphi\|_{L^\infty} )
        \label{lkslkcdhkhsdedj}
    \end{align}
    where $\chi_1$ and $\chi_2$ are bounded functions converging towards zero as $\epsilon\leq H/(C\ln \frac{1}{H})$,
    $\frac{h}{\epsilon}\downarrow 0$,  $\frac{\epsilon}{h} H \downarrow 0$ and $(\frac{h}{\epsilon})^2 \frac{1}{H} \downarrow 0$ (and $\Delta t\downarrow 0$ for $\chi_2$).
\end{Theorem}
\paragraph{Recall notations:}
$NK=L$ is the fixed spatial width, $[\eta^t]^x$ and $[\eta^t]^{-1}$ respectively denote the $x$ (slow) component and the inverse of the diffeomorphism $\eta^t$ (defined in Condition \ref{ContinuousFLAVORcondition}), $x(t)=[\eta^t]^x(u_1(t),\ldots,u_N(t))$ corresponds to the slow component of the exact PDE solution sampled at grid points.  $u_i(t)$ and $\bar{u}_i(t)$ represent the exact and the FLAVOR approximation of the solution to the semi-discrete system with the remainders \eqref{semi-discrete_exact}.

\begin{proof}
    The proof of Theorem \ref{ODEFLAVORTheorem} is analogous to that of Theorem 1.2 of \cite{FLAVOR10} (which will not be repeated here). The proof requires \eqref{TemporalConsistencyExact}, which is guarantied from Condition \ref{ConditionTemporalConsistency} by Lemma \ref{LemmaConsistentTemporal}.
    It is easy to check that the slow dependence on time of $f$, $g$, $\eta$ and $\mu$ does not affect the proof given in \cite{FLAVOR10}.
\end{proof}
\begin{Remark}
Condition \ref{ConditionTemporalConsistency} implies that the constant $C$ in Theorem \ref{ODEFLAVORTheorem} does not depend on $N$ or $K$.
This is important because although using a finer mesh leads to a smaller $K$ and a larger $N=L/K$, Condition \ref{ConditionTemporalConsistency} (which is equivalent to the accuracy of the semi-discrete approximation of the PDE) ensures that, as long as $K=\mathcal{O}(H)$ and $h \gg \epsilon H$, the constant $C$  in the error bounds on the slow component \eqref{jkdsgkdjshgdsdkjghe} and the fast component \eqref{lkslkcdhkhsdedj} will not blow up.
\end{Remark}

\begin{Remark}
  Observe that the application of the FLAVOR strategy  does not require the identification of the diffeomorphism $\eta$ (which may depend on the spatial discretization).
\end{Remark}

\section{On FLAVORizing characteristics}\label{seccharact}
The convergence result of the previous section is based on the semi-discretization of the original PDE.
PDEs and ODEs are also naturally connected via the method of characteristics, and henceforth it is natural to wonder whether a numerical integration of those characteristics by FLAVORs would lead to an accurate approximation of the solution of the original PDE. The answer to this question will be illustrated by analyzing  the following (generic)  PDE:
\begin{equation}
\begin{cases}\label{eqPDEcharac}
    F(Du,u,q,\epsilon^{-1})=0, & q\in U \\
    u(q)=\gamma(q), &q\in\Gamma
\end{cases}
\end{equation}
where $U\subset \mathbb{R}^d$ is the domain in which solution is defined,  $\Gamma$ and $\gamma$ define initial/boundary conditions.

The following condition corresponds to assuming that characteristics are well-posed.
\begin{Condition} Assume that
\begin{enumerate}
\item The PDE $F(Du,u,q,\epsilon^{-1})=0$ admits characteristics:
    \begin{eqnarray}
        \dot{q}&=&f(q,z,\epsilon^{-1}) \label{characteristics}  \\
        \dot{z}&=&g(q,z) \label{characteristics2} \\
        u(q(t))&=&z(t)
    \end{eqnarray}
    where $q\in U$ is a vector corresponding to coordinates of characteristics in the domain of the PDE, and $z$ corresponds to the unknown's value along the characteristics.
    \label{characteristicsConditionItem1}
\item
    For arbitrary $\epsilon$, any point in $U$ is reachable from the initial condition via one and only one characteristics.
\end{enumerate}
\label{characteristicsCondition}
\end{Condition}
The following conditions correspond to the assumption of existence of (possibly hidden) slow and locally ergodic fast variables for those characteristics.
\begin{Condition}
    Consider ODE \eqref{characteristics}. Assume that:
    \begin{enumerate}
    \item
         There exists a $z$-dependent diffeomorphism $\eta^z:q \mapsto [x,y]$ from $\mathbb{R}^d$ onto $\mathbb{R}^{d-p}\times \mathbb{R}^p$ with uniformly bounded $C^1, C^2$ derivatives with respect to both $q$ and $t$, such that $(x,y)$ satisfies (with $z(t)$ given by \eqref{characteristics2})
        \begin{equation}
        \begin{cases}
            \dot{x}&=f_1(x,y,z) \\
            \dot{y}&=\epsilon^{-1} f_2(x,y,z)
        \end{cases}
        \end{equation}
        where $f_1$, $f_2$, and $g$ have bounded $C^1$ derivatives with respect to $x$, $y$ and $z$, and $u([\eta^z]^{-1}(x,y))$ has bounded $C^1$ derivatives with respect to the (slow) variables $x$ and $z$.
    \item
         There exists a family of probability measures $\mu^z(x,dy)$ on $\mathbb{R}^p$ indexed by $x\in\mathbb{R}^{d-p}$ and $z\in\mathbb{R}$, as well as a family of positive functions $T \mapsto E^z(T)$ satisfying $\lim_{T\rightarrow \infty}E^z(T)=0$, such that for all $x_0,y_0,z_0,T$ bounded and $\phi$ uniformly bounded and Lipschitz, the solution to \begin{equation}
            \dot{Y}_t=f_2(x_0,Y_t,z_0)\quad \quad Y_0=y_0
        \end{equation}
        satisfies
        \begin{equation}
            \Big|\frac{1}{T}\int_0^T \phi(Y_s)ds-\int_{\mathbb{R}^p} \phi(y)\mu^{z_0}(x_0,dy)\Big|\leq \chi^{z_0}\big(\|(x_0,y_0)\|\big) E^{z_0}(T) (\|\phi\|_{L^\infty}+\|\nabla \phi\|_{L^\infty} )
        \end{equation}
        where $r \mapsto \chi^{z_0}(r)$ is bounded on compact sets, and $\mu^z$ has bounded derivative with respect to $z$ in total variation norm.
    \end{enumerate}
    \label{ContinuousFLAVORcondition2}
\end{Condition}
The second item of Condition \ref{ContinuousFLAVORcondition2} corresponds to the assumption that
the fast variable $y$ is locally ergodic with respect to a family of measures $\mu$ drifted by the slow variables $x$ and $z$.

The following lemma shows that, under the above conditions, the solution of PDE \eqref{eqPDEcharac} is nearly constant on the  orbit  of
the fast components ($y$) of the characteristics.

\begin{Lemma}
    Under Conditions \ref{characteristicsCondition} and \ref{ContinuousFLAVORcondition2}, for any fixed constant $C_1$ (independent of $\epsilon^{-1}$), there exists a constant $C_2$ independent of $\epsilon^{-1}$, such that for any $0 \leq t_1\leq C_1$, $0 \leq t_2\leq C_1$ and
    (fixed) $x_0$ and $z_0$,
    \begin{equation}
        \left| u\left([\eta^{z_0}]^{-1}(x_0,Y(t_1))\right)-u\left([\eta^{z_0}]^{-1}(x_0,Y(t_2))\right) \right| \leq C_2\epsilon
    \end{equation}
    where $Y(t_1)$ and $Y(t_2)$ are two points on the orbit of  $\dot{Y}(t)= f_2(x_0,Y(t),z_0)$.
    \label{OrbitalInvariance}
\end{Lemma}
\begin{proof}
    Under Conditions \ref{characteristicsCondition} and \ref{ContinuousFLAVORcondition2}, it is known (we refer for instance to \cite{MR810620} or to Theorem 14, Section 3 of Chapter II of \cite{MR1020057} or to \cite{MR2382139}) that $x$ and $z$ converge as $\epsilon\rightarrow 0$ towards $\tilde{x}$ and $\tilde{z}$ defined as the solution to the following ODEs with initial condition $x_0$ and $z_0$
    \begin{equation}
    \begin{cases}
        \dot{\tilde{x}}&=\int f_1 (\tilde{x},y,\tilde{z}) \mu^{\tilde{z}}(\tilde{x},dy) \\
        \dot{\tilde{z}}&=\int g ([\eta^{\tilde{z}}]^{-1}(\tilde{x},y),\tilde{z}) \mu^{\tilde{z}}(\tilde{x},dy)
        \label{averagedCharacteristics}
    \end{cases}
    \end{equation}
    Therefore, writing $y(t)$ the solution of $\dot{y}=\epsilon^{-1} f_2(\tilde{x},y,\tilde{z})$, we have as $\epsilon \rightarrow 0$
    \begin{equation}
        u([\eta^{\tilde{z}(t)}]^{-1}(\tilde{x}(t),y(t))) \rightarrow \tilde{z}(t)
    \end{equation}
    Now, taking the time derivative of $\hat{u}=u \circ \eta^{-1}$, we obtain
    \begin{equation}
        \hat{u}_x \dot{\tilde{x}} + \hat{u}_y \dot{y} + \hat{u}_z \dot{\tilde{z}} = \dot{\tilde{z}} + \dot{R}(\epsilon)
        \label{adsfhasglhglihrt}
    \end{equation}
    where $R(\epsilon)$ is a function of $t$ that goes to $0$ as $\epsilon\rightarrow 0$.

Furthermore,
    \begin{align*}
        \dot{Y}(t)&=f_2(x_0,Y(t),z_0) \\
        &=f_2(\tilde{x}(\epsilon t),y(\epsilon t),\tilde{z}(\epsilon t))+\frac{\partial f_2}{\partial \tilde{x}} (\tilde{x}(\epsilon t)-x_0)+\frac{\partial f_2}{\partial \tilde{z}} (\tilde{z}(\epsilon t)-z_0)+\frac{\partial f_2}{\partial y} (y(\epsilon t)-Y(t)) \\
        & \qquad \qquad +o(\epsilon)  + o(y(\epsilon t)-Y(t))
    \end{align*}
    By Taylor expansion, $\tilde{x}(\epsilon t)-x_0$ and $\tilde{z}(\epsilon t)-z_0$ are obviously $\mathcal{O}(\epsilon)$.
     Applying Gronwall's lemma, we also obtain that $y(\epsilon t)-Y(t)=\mathcal{O}(\epsilon)$. Therefore,
    \begin{equation}
        \dot{Y}(t)=f_2(\tilde{x}(\epsilon t),y(\epsilon t),\tilde{z}(\epsilon t))+\mathcal{O}(\epsilon)=\epsilon \dot{y}(t) + o(\epsilon)
        \label{g8orhuflbqhefqlbewfuhq}
    \end{equation}
Combining Eq. \ref{adsfhasglhglihrt} with Eq. \ref{g8orhuflbqhefqlbewfuhq}, we obtain
    \begin{align}
        & u\left(\eta^{-1}(x_0,Y(t_1))\right)-u\left(\eta^{-1}(x_0,Y(t_2))\right) = \int_{t_1}^{t_2} \hat{u}_y \cdot \dot{Y}(t) \, dt = \epsilon \int_{t_1}^{t_2} \hat{u}_y \cdot \dot{y} \, dt + o(\epsilon) \nonumber\\
        & \qquad = \epsilon \left( \int_{t_1}^{t_2} (\dot{\tilde{z}} - \hat{u}_x \dot{\tilde{x}} - \hat{u}_z \dot{\tilde{z}}) \, dt + R(\epsilon) \Big|^{t_2}_{t_1} \right) + o(\epsilon)
    \end{align}
    Since $\hat{u}_x$, $\dot{\tilde{x}}$, $\hat{u}_t$ and $\dot{\tilde{z}}$ are bounded, and $R(\epsilon)$ is vanishing (and hence bounded), we conclude that the right hand side is $\mathcal{O}(\epsilon)$.
\end{proof}

\begin{Condition}
    Assume that the domain $U$ is bounded (independently from $\epsilon^{-1}$).
    \label{boundedDomainCondition}
\end{Condition}

\begin{Lemma}
    If Conditions \ref{characteristicsCondition}, \ref{ContinuousFLAVORcondition2}, and \ref{boundedDomainCondition} hold, then every point in $U$ is reachable by a characteristics from the initial condition in bounded time  (independently from $\epsilon^{-1}$).
    \label{finiteTimeLemma}
\end{Lemma}

\begin{proof}
    From Condition \ref{characteristicsCondition}, we  already know that every point is reachable, and therefore it suffices to show
    that hitting times do not blow up as $\epsilon\rightarrow 0$. Since $x(\cdot)$ converges to $\tilde{x}(\cdot)$ (see proof of Lemma \ref{OrbitalInvariance}), by considering the $x$ component of the characteristics (projected by $\eta$), it becomes trivial to show that the hitting time converges to a fixed value (and hence,  does not blow up).
     Using Condition \ref{boundedDomainCondition}, we conclude that that any point in
      $U$ can be hit in (uniformly) bounded time from the initial condition.
\end{proof}

 Analogously to the Integrator \ref{FLAVOR_ODE_integrator}, a legacy integrator for \eqref{characteristics} and \eqref{characteristics2} can be FLAVORized, and shown to be convergent under regularity and stability conditions (analogous to Condition \ref{RegularityCondition}) requiring $f_1$, $f_2$ and $g$ to be Lipschitz continuous and $\tilde{q}(t)$ and $\tilde{z}(t)$ to be bounded. The convergence result is analogous to Theorem  \ref{ODEFLAVORTheorem}, modulo the following change of notation: the slow index is now $z$ instead of $t$, the original coordinates are $q$ instead of $u_i$, the vector field of the original coordinates is $f$ instead of $f_i$, and the dynamics of the slow index comes from the non-trivial drift of $\dot{z}=g(q,z)$ instead of the trivial $\dot{t}=1$. We define $\tilde{u}(\tilde{q}(t)):=\tilde{z}(t)$ for all $t$ on each FLAVORized characteristics $[\tilde{q}(t),\tilde{z}(t)]$. Naturally, $\tilde{u}$ is only defined at discrete points in the domain $U$. These discrete points, however, densely `fill' the space in the sense that (as shown by the  proof of the following theorem) FLAVORied characteristics remain very close to exact characteristics ($x$ components are close in Euclidean distance, and $y$ components are close as well in terms of orbital distance induced by the infimum of point-wise Euclidean distances).

By the two-scale convergence theorem, we can quantify: the strong convergence of the slow coordinate of the characteristics and the unknown's value along the characteristics, and  the weak convergence of fast coordinate of the characteristics. Finally, these single characteristics' ODE approximation error bounds can be transferred to the PDE approximation error bounds by considering the entire family of characteristics starting from all points (in initial condition).

\begin{Theorem}
Write $\tilde{u}(\tilde{q})$ the solution obtained by FLAVORizing all characteristics.
    Under Conditions \ref{characteristicsCondition}, \ref{ContinuousFLAVORcondition2}, \ref{boundedDomainCondition},  the consistency and regularity and stability Conditions corresponding to Conditions \ref{ConditionTemporalConsistency}  and \ref{RegularityCondition} (with the change of notation described above),  there exist a constant $C$ independent of $\epsilon^{-1}$ and $q_0\in \Gamma$, such that \begin{equation}
        |\tilde{u}(\tilde{q})-u(\tilde{q})|\leq C \chi_1(q_0,\gamma(q_0),\epsilon,\delta,\tau) (1+\chi_2(q_0,\gamma(q_0),\epsilon,\delta,\tau,T,t))
    \end{equation}
    for any $\tilde{q}$ on any FLAVORized characteristics, where $q_0\in \Gamma$ and $\gamma(q_0)$ correspond to the initial condition that leads to $\tilde{q}$ via a FLAVORized characteristics, and $\chi_1$ and $\chi_2$ are vanishing error bound functions.
\end{Theorem}
\begin{Remark}
When $\Gamma$ is compact (such as in the case of periodic boundary condition), $\chi_1$ and $\chi_2$ can be further chosen to be independent of $q_0$ (hence $\tilde{q}$) by taking a supremum  over $\Gamma$.
\end{Remark}

\begin{proof}
    By Condition \ref{characteristicsCondition}, all $q\in U$ can be traced back to $q_0\in\Gamma$ through a characteristics. By Lemma \ref{finiteTimeLemma},   characteristics starting from $q_0$ reach $q$ in bounded time $T$. Using the two-scale convergence of the FLAVORization of these characteristics (a result analogous to Theorem \ref{ODEFLAVORTheorem}), we deduce that the approximation error associated with $\tilde{z}_T$ (on FLAVORized characteristics) can be bounded $C\chi_1$ (with respect to the true value $u(q)=z_T$, the error $Ce^{CT}$ has been replaced by $C$ because $T$ is bounded).

    Now observe  that $\tilde{q}_T \neq q_T$, where $\tilde{q}_T$ is the coordinate of the FLAVORized characteristics starting from $q_0$. As before, let $[x_T,y_T]=\eta(q_T)$ and $[\tilde{x}_T,\tilde{y}_T]=\eta(\tilde{q}_T)$. The error on the slow component is $\| x_T-\tilde{x}_T\| \leq C \chi_1$. The possible large error on the fast component is not a problem because we can look for a near-by point on the fast orbit with introducing only an $\mathcal{O}(\epsilon)$ error on the unknown's value (Lemma \ref{OrbitalInvariance}):
    \begin{equation}
    \begin{cases}
        u(\eta(x_T,y_T))=u(\eta(x_T,y^*_T))+\mathcal{O}(\epsilon) \\
        y^*_T=\arg \min_{Y_t | \dot{Y}_t=f(x_T,Y_t)} \|\tilde{y}_T-Y_t \|
    \end{cases}
    \end{equation}
    Since $\|\tilde{x}_T-x_t\|$ is small, the local ergodic measures that represent the orbits given by $\dot{Y}_t=f(x_T,Y_t)$ and $\dot{Y}_t=f(\tilde{x}_T,Y_t)$ will be small: $\|\mu(x_T,dy)-\mu(\tilde{x}_T,dy)\|_{\text{T.V.}} \leq C \chi_1 \chi_2$ is by chain rule. Because $\tilde{y}_T$ is on the orbit of $\dot{Y}_t=f(\tilde{x}_T,Y_t)$, we will have $\|y^*_T-\eta^y(\tilde{q}_T)\|\leq C \chi_1 \chi_2$.

     All together, we obtain
     \begin{align}
        |\tilde{u}(\tilde{q}_T)-u(\tilde{q}_T)| &= |\tilde{z}_T-u(\tilde{q}_T)| \nonumber\\
        &\leq |\tilde{z}_T-u(q)| + |u(q_T)-u(\tilde{q}_T)| \nonumber\\
        &\leq C\chi_1 + C \| \nabla (u\circ \eta) \|_\infty \left( \|x_T-\eta^x(\tilde{q}_T)\|+\|y_T-\eta^y(\tilde{q}_T)\| \right) \nonumber\\
        &\leq C \chi_1 + C (\chi_1 + \chi_1 \chi_2) = C \chi_1 + C \chi_1 \chi_2
     \end{align}
\end{proof}

\begin{Remark}
To keep the presentation concise, we have written $C$ all constants that do not depend on essential parameters.
\end{Remark}

\begin{Remark}
    As shown above, $u$ will be captured strongly. $Du$, on the other hand, depends on a derivative with respect to the fast variable, and therefore will only be convergent in a weak sense.
\end{Remark}

\paragraph{Relevance to an error analysis for PDE-FLAVORS}
The above result guarantees the convergence of FLAVORized characteristics. It is also possible to establish an error bound on the difference between a specific PDE-FLAVOR discretization and the approximation given by the above FLAVORized characteristics (and hence prove the convergence of this specific PDE-FLAVOR discretization). Such an error bound  could be obtained by first transforming FLAVORized characteristics to PDE-FLAVOR grid points via interpolating functions, and then using the fact that  coordinate transformations do not affect the efficiency of FLAVORS. For the sake of conciseness, we did not elaborate on this point here.

\section{Acknowledgement}
This work is supported by NSF grant CMMI-092600. We  thank Guo Luo for stimulating discussions and Sydney Garstang for proofreading the manuscript.

\bibliographystyle{siam}
\bibliography{molei15}

\end{document}